\documentclass[11pt,reqno]{amsart}
\usepackage{amsmath,amsthm,amsfonts,amssymb}
\usepackage{mathtools}
\usepackage[czech,english]{babel}
\usepackage{graphicx}
\usepackage[all]{xy}
\usepackage{longtable}
\usepackage[format=hang,list=off]{caption}
\usepackage{tensor}
\usepackage[shortlabels]{enumitem}
\usepackage{constants}
\usepackage{color}
\definecolor{shadecolor}{gray}{.85}%
\definecolor{tintedcolor}{gray}{.80}%
\usepackage{framed}%
\definecolor{mytintedcolor}{gray}{.95}%
\makeatletter
\newdimen\svparindent

\newenvironment{mytinted}{%
  \MakeFramed {\FrameRestore}}%
{\endMakeFramed}
{\endlist\end{mytinted}\egroup}
\makeatother
\usepackage[OT1]{fontenc}
\newtheorem{theorem}{Theorem}

\newtheorem{example}{Example}
\newtheorem{lemma}{Lemma}

\newtheorem{corollary}{Corollary}

\newtheorem{claim}{Claim}
\newtheorem{definition}{Definition}

\newtheorem{remark}{Remark}
\newcommand{\bbbn}{\mathbb{N}}

\newcommand{\bbbr}{{\mathbb R}}

\newcommand{\Type}[2]{{\sf\textstyle Type}^{\scriptscriptstyle #1}_{#2}}
\newcommand{\restr}[2]{{#1}\upharpoonright_{#2}}
\newconstantfamily{r}{symbol=r}
\newconstantfamily{eps}{symbol=\epsilon}
\newconstantfamily{N}{symbol=N}
\DeclareMathOperator\qrank{\rm qrank}
\DeclareMathOperator\lrank{\rm lrank}
\newcommand\dist[2]{{\rm Dist}_{#1,#2}}
\newcommand\ldist[2]{{\rm Dist}_{#1,#2}^{\rm local}}

\title{Approximations of Mappings}
\thanks{Supported by grant ERCCZ LL-1201 
and CE-ITI, and by the European Associated Laboratory ``Structures in
Combinatorics'' (LEA STRUCO) P202/12/G061}
\author[J. Ne{\v s}et{\v r}il]{Jaroslav Ne{\v s}et{\v r}il}
\address{Jaroslav Ne{\v s}et{\v r}il\\
Computer Science Institute of Charles University (IUUK and ITI)\\
   Malostransk\' e n\' am.25, 11800 Praha 1, Czech Republic}
\email{nesetril@kam.ms.mff.cuni.cz}
\author[P. Ossona de Mendez]{Patrice~Ossona~de~Mendez}
\address{Patrice~Ossona~de~Mendez\\
Centre d'Analyse et de Math\'ematiques Sociales (CNRS, UMR 8557)\\
  190-198 avenue de France, 75013 Paris, France
  and
     Computer Science Institute of Charles University (IUUK)\\
   Malostransk\' e n\' am.25, 11800 Praha 1, Czech Republic}
 \email{pom@ehess.fr}
 \date{\today}
 


\begin{document}
 \begin{abstract}
 We consider mappings, which are structure consisting of a single function (and possibly some number of unary relations) and address  the problem of approximating a continuous mapping by a finite mapping. This problem is the inverse problem of the construction of a continuous limit for first-order convergent sequences of finite mappings. We solve the approximation problem and, consequently, the full characterization of limit objects for mappings for first-order (i.e. ${\rm FO}$) convergence and local (i.e. ${\rm FO}^{\rm local}$) convergence. 
 
 This work can be seen both as a first step in the resolution of inverse problems (like Aldous-Lyons conjecture) and a strengthening  of the classical  decidability result for finite satisfiability in Rabin class (which consists of first-order logic with equality, one unary function, and an arbitrary number of monadic predicates).
 
 The proof involves model theory and analytic techniques.
 
 \end{abstract}
\maketitle

\renewcommand{\contentsname}{\vspace{-2\baselineskip}}
\setcounter{tocdepth}{1}
\tableofcontents

\section{Introduction}
\label{sec:intro}
We consider the following {\em approximation problems}:
Given an infinite structure with given first-order properties, as well as satisfaction probabilities for every first-order formula, can one find a finite structure with approximately similar properties and satisfaction probabilities? What if we are not given the infinite structure, but only the satisfaction probability of first-order formulas? 

These problems are in general intractable, as (even when considering no probabilities of satisfaction) it is known that deciding whether a sentence satisfied by an infinite structure is also satisfied by a finite structure is (in general) undecidable. Intensive studies have been conducted to determine decidable classes of structures and fragments of first-order logic. A maximal example is the {\em Rabin class}, which consists of all first-order sentences with arbitrary quantifier prefix and equality, one unary function symbol, and an arbitrary number of unary relation symbols (but no function or relation symbols of higher arity). The satisfiability problem and the finite satisfiability problem for this class are both decidable, but not elementary recursive \cite{borger2001classical}.

Another particular case of our problem was considered extensively in the context of topological group theory, ergodic theory and graph limits, and concerns the class of bounded degree graphs (one binary symmetric symbol) and local first-order formulas with a single free variable.  It can be formulated as follows: consider a unimodular probability measure $\mu$ defined on the set $\mathcal{G}^*$ of all countable rooted connected graphs endowed with the metric defined by the rooted neighborhood isomorphisms. Can $\mu$ be approximated by finite graphs? This question is known as the {\em Aldous--Lyons conjecture}. It is not just an isolated problem as a positive solution would have far-reaching consequences, by proving that all finitely generated groups are sofic (answering a question by Weiss \cite{Weiss2000}),  the direct finiteness conjecture of Kaplansky \cite{kaplansky1972fields} on group algebras, a conjecture of Gottschalk \cite{Gottschalk1973} on surjunctive groups in topological dynamics, the Determinant Conjecture on Fuglede-Kadison determinants, and Connes' Embedding Conjecture for group von Neumann algebras \cite{Connes1976}.
It is easily shown that Aldous-Lyons conjecture can be reduced to the approximation problem for quantifier-free formulas on structures with two functions $f$ and $g$ satisfying $f^2=g^3={\rm Id}$. 

In this paper we solve the approximation problem for mappings, i.e. structures consisting of a set $X$ and an (endo)function $f : X \to X$, and more generally we solve it for the whole Rabin class.  
 At lest at first glance it is perhaps surprising that 
 such a seemingly special case is quite  difficult to handle.

Approximation problems recently appeared in the context of graph limits as so called {\em inverse problems}.
In order to make the connection 
clear, we take time for a quick review of some of the fundamental notions and problems encountered in the domain of graph limits, and how they are related to the study of limits and approximations of algebras (that is of functional structures). 

A sequence of (colored) graphs with maximum degree at most $d$ converges if, for every integer $r$, the distribution of the isomorphism type of the ball of radius $r$ rooted at a random vertex (drawn uniformly at random) converges. The limit object of a local convergent sequence of graphs is a {\em graphing}, that is a graph on a standard Borel space, which satisfies a {\em  Mass Transport Principle}, which amounts to say that for every Borel subsets $A,B$ it holds that
$$
\int_A{\rm deg}_B(v)\,{\rm d}v=\int_B{\rm deg}_A(v)\,{\rm d}v.
$$
	An alternative description of graphings is as follows: a graphing is defined by a finite number of measure preserving involutions $f_1,\dots,f_D$ on a standard Borel space, which define the edges of the graphing as the union of the orbits of size two of 
$f_1,\dots,f_D$. 

The idea to conceptualize limits of structures by means of convergence of the satisfaction probability of formulas in a fixed fragment of first-order logic has been introduced by the authors in 
\cite{CMUC}.
In this setting, a sequence $(\mathbf A_n)_{n\in\bbbn}$ of structures is {\em convergent} (or {\em $X$-convergent}) if, for every first-order formula $\phi$ in a fixed fragment $X$ the probability $\langle\phi,\mathbf A_n\rangle$ that $\phi$ is satisfied in $\mathbf A_n$ for a random assignment of elements of $\mathbf A_n$ to the free variables of $\phi$ converges as $n$ grows to infinity. If $X$ is the set of all first-order formulas, then we speak about FO-convergence.
This definition allowed us to consider limits of general combinatorial structures, and was applied to limits of sparse graphs with unbounded degrees \cite{limit1,modeling,Gajarsky2015,modeling_jsl}, 
  matroids \cite{matroid_limit}, and tree semi lattices \cite{QFTSL}. 

The  main result of \cite{MapLim} is the construction of a limit object for ${\rm FO}$-convergent sequences of mappings (a {\em mapping} being an algebra with a single function symbol and --- possibly --- finitely many unary predicates).
\begin{theorem}
\label{thm:folimmap}
Every  ${\rm FO}$-convergent sequence $(\mathbf F_n)_{n\in\bbbn}$ of finite mappings  (with $\lim_{n\rightarrow\infty} |F_n|=\infty$) has a modeling mapping limit $\mathbf L$, such that
\begin{enumerate}
	\item the probability measure $\nu_{\mathbf L}$ is atomless;
	\item the complete theory of $\mathbf L$ has the finite model property;
	\item $\mathbf L$ satisfies the finitary mass transport principle.
\end{enumerate} 
\end{theorem}

Let us explain the (undefined) notions appearing in this theorem:
\begin{enumerate}[(i)]
\item A {\em modeling} $\mathbf L$ is a {\em totally Borel structure} --- that is a structure whose domain $L$ is a standard Borel space, such that every definable set is Borel --- endowed with a probability measure $\nu_{\mathbf L}$.
	\item
	 The measure $\nu_{\mathbf L}$ is {\em atomless}
	(or {\em continuous}, or {\em diffuse})  if for every $v\in L$ it holds $\nu_{\mathbf L}(\{v\})=0$.
	(As we consider only standard Borel spaces, this condition is equivalent to the condition that for every Borel subset $A$ with $\nu_{\mathbf L}(A)>0$ there exists a Borel subset $B$ of $A$ with $\nu_{\mathbf L}(A)>\nu_{\mathbf L}(B)>0$.)
	 The necessity of this condition is witnessed by the formula $x_1=x_2$, as $\langle x_1=x_2,\mathbf F\rangle=1/|F|$ holds
for every finite mapping $\mathbf F$. This conditions is thus required as soon as we consider ${\rm QF}$-convergence.
	\item the {\em finitary mass transport principle} (FMTP) means that for every Borel subsets $X,Y$ of $L$ and every positive integer $k$  it holds 
\begin{align*}
	(\forall v\in Y)\,|f^{-1}(v)\cap X|=k\quad&\Rightarrow\quad \nu_{\mathbf L}(f^{-1}(Y)\cap X)=k\nu_{\mathbf L}(Y)\\
	(\forall v\in Y)\,|f^{-1}(v)\cap X|>k\quad&\Rightarrow\quad \nu_{\mathbf L}(f^{-1}(Y)\cap X)>k\nu_{\mathbf L}(Y)\\
\end{align*}
This condition can be reformulated as follows: 
the set of all $y$ such that $f_{\mathbf F}^{-1}(y)$ is infinite has zero $\nu_{\mathbf F}$-measure, and 
for every Borel subsets $X$ and $Y$ of $L$ (with $|f_{\mathbf F}^{-1}(y)|<\infty$ for all $y\in Y$) we have 
\begin{equation}
\label{eq:fmtp}
\nu_{\mathbf F}(X\cap f_{\mathbf F}^{-1}(Y))=
\int_Y |f_{\mathbf F}^{-1}(y)\cap X| \,{\rm d}\nu_{\mathbf F}(y).	
\end{equation}

When $X$ and $Y$ are definable subsets, the above condition is clearly required for being a limit.
	\item the {\em finite model property} 
	means that for every sentence $\theta$ satisfied by $\mathbf L$ there exists a finite mapping $\mathbf F$ that satisfies $\theta$. This is indeed a necessary condition for $\mathbf L$ to be an elementary limit of finite mappings hence necessary as soon as we consider
	${\rm FO}$-convergence. As mentioned, the problem of existence of a  finite mapping $\mathbf F$  satisfying a given sentence $\theta$ is decidable, though with huge time complexity.
\end{enumerate}

Theorem~\ref{thm:folimmap} was proved as a combination of general results about limit distributions from \cite{CMUC} 
and methods developed in \cite{modeling} for the purpose of graph-trees. This theorem
has the following corollary.
\begin{corollary}
\label{cor:local}
Every  ${\rm FO}^{\rm local}$-convergent sequence $(\mathbf F_n)_{n\in\bbbn}$ of finite mappings  (with $\lim_{n\rightarrow\infty} |F_n|=\infty$) has a modeling mapping ${\rm FO}^{\rm local}$-limit $\mathbf L$, such that
\begin{enumerate}
	\item the probability measure $\nu_{\mathbf L}$ is atomless;
	\item $\mathbf L$ satisfies the finitary mass transport principle.
\end{enumerate} 	
\end{corollary}
\begin{proof}
	Consider an ${\rm FO}$-convergent subsequence. Such a subsequence exists by (sequential) compactness of ${\rm FO}$-convergence. According to Theorem~\ref{thm:folimmap} this subsequence has a modeling mapping limit $\mathbf L$ satisfying all the requirements. This modeling limit is then a modeling ${\rm FO}^{\rm local}$-limit of $(\mathbf F_n)_{n\in\bbbn}$.
\end{proof}

The inverse problems aim to determine which objects are limits of finite mappings (for given types of convergence). 
The main contribution of this paper is the solution of the inverse problems for mappings. Namely, for FO and ${\rm FO}^{\rm local}$-convergence we show how to approximate a modeling mapping by a finite mapping. (For QF-convergence the inverse problem is much easier and was solved in \cite{MapLim}.)

\begin{theorem}
\label{thm:invlocmap}
Every atomless modeling mapping $\mathbf L$  that satisfies the finitary mass transport principle is the 
${\rm FO}^{\rm local}$-limit of an ${\rm FO}^{\rm local}$-convergent sequence of finite mappings.
\end{theorem}

\begin{theorem}
\label{thm:invfomap}
Every atomless modeling mapping $\mathbf L$ with the finite model property that satisfies the finitary mass transport principle is the 
${\rm FO}$-limit of an ${\rm FO}$-convergent sequence of finite mappings.
\end{theorem}

Here is a rough outline of the proof of
Theorem~\ref{thm:invfomap}:
\begin{enumerate}
	\item reduce to the case where the mapping modeling $\mathbf L$ has no connected component of measure greater than $\epsilon$;
	\item consider a derived modeling mapping $\mathbf L'$ obtained by removing all the elements with zero-measure rank-$R$ local type;
	\item cut all the short circuits by means of interpretation;
	\item approximate the measure on the rank-$R$ local types by a rational measure $\mu$;
	\item construct a finite mapping  $\mathbf F$ such that the measure of each rank-$r$ local type $t$ is equal to 
	what is derived from $\mu$;
	\item\label{step:merge} consider a finite mapping $\mathbf M$, which is equivalent to $\mathbf L$ up to a huge quantifier rank, and
	 merge it with a great number of copies of $F$ to form an ${\rm FO}$-approximation of $\mathbf L$;
	\item deduce, using interpretation, an ${\rm FO}$-approximation of the original mapping modeling.
\end{enumerate}

Theorem~\ref{thm:invlocmap} is then proved by considering separately large and small connected components, and following a similar strategy as the proof of Theorem~\ref{thm:invfomap}:
\begin{enumerate}
	\item every connected modeling mapping is close (in the sense of local convergence) to a modeling mapping with finite height; such a modeling mapping has the finite model property, hence maybe ${\rm FO}$-approximated thanks to Theorem~\ref{thm:invfomap};
	\item[(2-5)]\addtocounter{enumi}{4} for a mapping modeling without connected components of measure greater than $\epsilon$, construct a finite mapping $\mathbf F$ as in the steps (2) to (5) of the proof of Theorem~\ref{thm:invfomap};
	
	\item then complete $\mathbf F$ by means of small models of missing necessary local types, merged with a great number of copies of $\mathbf F$;
	\item the ${\rm FO}^{\rm local}$ approximation is obtained as the disjoint union of the ${\rm FO}^{\rm local}$ approximations of large connected components and the ${\rm FO}_1^{\rm local}$ approximation of the remaining components (after careful tuning of the respective orders).
\end{enumerate}  

It should be noted that Theorems~\ref{thm:invlocmap} and~\ref{thm:invfomap} allow to obtain approximations from a mapping modeling, which may have only finitely many unary predicates in its signature. The case where we allow infinitely many unary predicates easily restricts to this case, as (for given metrization of ${\rm FO}$- and ${\rm FO}^{\rm local}$-convergence) for every $\epsilon>0$ there exist $\epsilon'>0$ and $C\in\mathbb N$ such that any $\epsilon'$-approximation of the mapping considering only the first $C$ unary predicates is an $\epsilon$-approximation of the mapping when considering all the unary predicates.
Hence Theorem~\ref{thm:invfomap} and~\ref{thm:invlocmap} solves the approximation problem for the Rabin class modelings.

As a pleasing consequence of our general methods we believe
 that Theorems~\ref{thm:invlocmap} and~\ref{thm:invfomap} can be formulated in a setting where we do not approximate a particular modeling $\mathbf L$ but rather consider the satisfaction probability of formulas. This would gives a full solution of the second type approximation problem for the Rabin class.

\section{Preliminaries}
\subsection{Facts from Finite Model Theory}
\label{sec:prelim_model}
We recall some basic definitions and facts from finite model theory. The interested reader is refereed to \cite{Ebbinghaus1996,Hodges1993,Hodges1997,Lascar2009,Libkin2004,Marker2001}.

A {\em signature} $\sigma$ is a list function or relation symbols with their arities. 
A {\em $\sigma$-structure} $\mathbf A$ is defined by its {\em domain}  $A$, its {\em signature} $\sigma$, and the interpretation in $A$ of all the relations and functions  in $\sigma$. The {\em Gaifman graph} of a $\sigma$-structure $\mathbf A$ is the graph with vertex set $A$, where two elements are adjacent if they belong to a same relation or are related by a function application. When we speak about the neighborhood of an element $x$ in $A$ or about the distance between two elements $x$ and $y$ in $A$, we mean the set of elements adjacent to $x$ in the Gaifman graph of $\mathbf A$ or the graph distance between $x$ and $y$ in the Gaifman graph of $\mathbf A$. Also, for $u\in A$ and $r\in\bbbn$ we denote by $B_r(\mathbf A,u)$ the {\em $r$-ball} of $u$ in $\mathbf A$, that is the set of all elements of $A$ at distance at most $r$ from $u$.

We denote by ${\rm FO}(\sigma)$ the set of all first-order formulas (in the language defined by the signature $\sigma$).
A formula $\phi$ (with $p$ free variables) is {\em local} if its satisfaction only depends on a fixed $r$-neighborhood of its free variables, and we denote by ${\rm FO}^{\rm local}(\sigma)$ the set of all local formulas.
Also, we denote by ${\rm QF}(\sigma)$ the set of all quantifier free formulas. When we consider sub-fragments where we restrict free variables to $x_1,\dots,x_p$, we will add $p$ as a subscript, as in 
${\rm FO}_p(\sigma)$ or ${\rm FO}_p^{\rm local}(\sigma)$.

For a first-order formula $\phi$ with $p$ free variables and a $\sigma$-structure $\mathbf A$ be define $\phi(\mathbf A)$ as the set of all $p$-tuples of elements of $\mathbf A$ that satisfy the formula $\phi$ in $\mathbf A$, that is:
$$\phi(\mathbf A)=\{(v_1,\dots,v_p)\in A^p:\ \mathbf A\models\phi(v_1,\dots,v_p)\}.$$

In the following definition we consider signatures with a function symbol and finitely many unary predicates. Although Rabin class allows infinitely many unary predicates, this is not a real restriction in the context of approximation problems, but this assumption will make the definitions and notations simpler.

\begin{definition}
\label{def:mapping}
A {\em mapping} is a $\sigma$-structure, where the signature  $\sigma$ consists of a single unary function symbol $f$ and (possibly) finitely many unary relation symbols $M_1,\dots,M_c$.
\end{definition}

Let $\mathbf F$ be a mapping. We denote by $F$ the domain of $\mathbf F$ and by $f_{\mathbf F}$ the interpretation of the symbol $f$ in $\mathbf F$ (thus $f_{\mathbf F}:F\rightarrow F$). Unary relations will be denoted by $M_i^{\mathbf F}$ (or simply just $M_i$).
Note that  the distance ${\rm dist}(u,v)$ between two elements $u,v$ in a mapping $\mathbf F$ is the minimum value $a+b$ such that $a,b\geq 0$ and $f_{\mathbf F}^a(u)=f_{\mathbf F}^b(v)$.

Every formula $\phi\in{\rm FO}_1^{\rm local}$ is logically equivalent to a formula with no function composition. Such formulas we call {\em clean}.	

\begin{definition}
\label{def:rank}
	The {\em quantifier rank}  of a formula $\phi$, denoted by  $\qrank(\phi)$, is the minimum number of nested quantifiers in a clean formula equivalent to $\phi$.
	
	The {\em local rank} of a local formula $\phi$, denoted by  $\lrank(\phi)$, is the minimum number of nested quantifiers in a clean formula equivalent to $\phi$ in which quantification is restricted to previously defined variables and their neighbors.
\end{definition}

It is easily checked that for a given finite signature $\sigma$ there exist only finitely many local formulas $\phi\in{\rm FO}_1^{\rm local}(\sigma)$ that have local rank at most $r$ (up to logical equivalence).

A {\em local type} is any maximal consistent subset $t$ of ${\rm FO}_1^{\rm local}(\sigma)$. 
	The {\em local type} of an element $v$ of a mapping $\mathbf F$ is the local type $t$ such that $\mathbf F\models \phi(v)$ holds for every $\phi\in t$. 
 A {\em rank $r$ local type} is the subset of a all formulas with rank at most $r$ in a local type. We denote by $\mathcal T_r(\sigma)$ the set of all rank $r$ local types for signature $\sigma$. 
 We denote by 
	$\Type{\mathbf F}{r}(v)$ the rank $r$ local type of an element $v$ in a mapping $\mathbf F$.
	
Note that for every rank $r$ local type $t\in  \mathcal T_r$ there exists a clean formula $\varphi_t\in t$ (in which quantification is restricted to previously defined variables and their neighbors) such that $\varphi_t$ is logically equivalent to the conjunction of all the formulas in $t$. (The formula $\varphi_t$ will always have this meaning.) Thus for every $\sigma$-structure $\mathbf F$ and every $v\in F$ it holds that
$$
\Type{\mathbf F}{r}(v)=t\quad\iff\quad \mathbf F\models\varphi_t(v).
$$

For $r<r'$, $t\in \mathcal T_r(\sigma)$ and $t'\in \mathcal T_{r'}(\sigma)$ we say that $t'$ {\em refines} $t$, and write $t'\prec t$,  if 
$\varphi_{t'}\vdash \varphi_t$ (i.e. if $t'\supseteq t$).

Given two mappings $\mathbf F$ and $\mathbf F'$, it is well known that $\mathbf F$ and $\mathbf F'$ satisfy the same sentences with quantifier rank at most $r$, what is denoted by $\mathbf F\equiv_r\mathbf F'$, if and only if Duplicator has a winning strategy for the $r$-rounds Ehrenfeucht--Fra{\"\i}ss\'e game. 

Given two elements $v\in F$ and $v'\in F'$, testing whether $\Type{\mathbf F}{r}(v)=\Type{\mathbf F'}{r}(v')$ can be done using a variant of a Ehrenfeucht--Fra{\"\i}ss\'e game: We start by defining $u_0=v$ and $u_0'=v'$. At each round $1\leq k\leq r$, Spoiler chooses in $F$ an element $u_k$ adjacent to some of $u_0,\dots,u_{k-1}$ (or
in $F'$ an element $u_k'$ adjacent to some of $u_0',\dots,u_{k-1}'$). Then Duplicator should choose $u_k'\in F'$ (or $u_k\in F$) so that for every $0\leq i,j\leq k$ it holds
\begin{align*}
\mathbf F\models u_i=u_j\quad&\iff\quad\mathbf F'\models u_i'=u_j'\\
\mathbf F\models f(u_i)=u_j\quad&\iff\quad\mathbf F'\models f(u_i')=u_j'
\end{align*}
	Spoiler wins if Duplicator cannot make such a choice and $k\leq r$; otherwise, Duplicator wins. It is easily checked that $\Type{\mathbf F}{r}(v)=\Type{\mathbf F'}{r}(v')$ if and only if Duplicator has a winning strategy. We call this variant of Ehrenfeucht--Fra{\"\i}ss\'e game the {\em local Ehrenfeucht--Fra{\"\i}ss\'e game}.
	
For $r \leq r'$ we define the natural projection $\pi_r$ mapping an $r'$-type $t$ to the $r$-type $\pi_r(t)$, which is just the subset of all formulas in $t$ with rank at most $r$. Obviously, if $r'>r$ then $\pi_r(\Type{\mathbf F}{r'}(v))=\Type{\mathbf F}{r}(v)$.

Let $\sigma,\sigma'$ be signatures of mappings.
Let $M_1,\dots,M_a$ be the symbols of the unary symbols in $\sigma'$
(as usual $f$ is the function symbol). The following is a standard definition.

\begin{definition}
\label{def:interp}
A {\em basic interpretation} $\mathsf I$ of $\sigma'$-structures into $\sigma$-structures is defined by $a$ formulas $\kappa_1,\dots,\kappa_a$ with a single free variable, and a formula $\eta$ with two free variables defining the graph of an endofunction, that is such that 
\[\vdash\ \forall x\ \exists y\  \bigl(\eta(x,y)\wedge (\forall z)(\eta(x,z)\rightarrow (z=y))\bigr).\]

For every $\sigma$-structure $\mathbf A$, the $\sigma'$-structure $\mathbf B=\mathsf I(\mathbf A)$ has same domain as $\mathbf A$ (i.e. $B=A$), its relations are defined by
$$\mathbf B\models M_i(v)\quad\iff\quad \mathbf A\models\kappa_i(v)$$
and $f_{\mathbf B}$ is (implicitly) defined by
$$\mathbf B\models f(u)=v\quad\iff\quad \mathbf A\models\eta(u,v).$$
The interpretation $\mathsf I$ is {\em trivial} if $\eta(x,y):=(f(x)=y)$ (hence $f_{\mathbf B}=f_{\mathbf A}$).
\end{definition}

For every first order formula $\phi$ with $p$ free variables (on the language of $\sigma'$-structures)  the first-order formula $\mathsf I(\phi)$ is obtained by replacing (in a clean formula logically equivalent to $\phi$) terms 
$M_i(x)$ by $\kappa_i(x)$ and terms $f(x)=y$ by $\eta(x,y)$.
The formula $\mathsf I(\phi)$ is  such that
for every $\sigma$-structure $\mathbf A$
and  every $v_1,\dots,v_p\in B$ it holds
$$
\mathbf B\models \phi(v_1,\dots,v_p)\quad\iff\quad\mathbf A\models\mathsf I(\phi)(v_1,\dots,v_p).
$$
Note that if $\phi$ and all the formulas defining $\mathsf I$ are local then $\mathsf I(\phi)$ is local and 
$${\rm lrank}(\mathsf I(\phi))\leq {\rm lrank}(\phi)+\max({\rm lrank}(\kappa_1),\dots,{\rm lrank}(\kappa_a),{\rm lrank}(\eta)).$$

\subsection{Structural Limits}
\label{sec:lim}
We recall here some definitions and notations from \cite{CMUC}.

Recall that a $\sigma$-structure is {\em Borel} if its domain is a standard Borel space, and all the relations and functions of the structure are Borel.
For instance, the mapping $\mathbf F$ is {\em Borel} if the function
$f_{\mathbf F}:F\rightarrow F$ and the subsets $M_i(\mathbf F)=\{v\in F: \mathbf F\models M_i(v)\}$ are Borel;

 A stronger notion has been proposed in \cite{CMUC}:

	\begin{definition}
	\label{def:modeling}
		A {\em $\sigma$-modeling} (or a {\em modeling} when $\sigma$ is implied) is a $\sigma$-structure $\mathbf M$, whose domain $M$ is a standard Borel space endowed with a 
		probability measure $\nu_{\mathbf M}$, and with the property that every definable subset of a power of $M$ is Borel.
	\end{definition}

If $\mathbf F$ is a finite structure, it will be practical to implicitly consider a uniform probability measure $\nu_{\mathbf F}$ on $F$, for the sake of simplifying the notations.

Note that every modeling mapping is obviously Borel, but the converse does not hold true in general, as shown by the next example.
\begin{example}
	\label{ex:souslin}
A counter-example of Lebesgue's belief that the projection to $\bbbr$ of a Borel subset of $\bbbr^2$ is Borel has been given by Souslin.
It follows that there exits a Borel subset $S\subseteq (0,1]\times (0,1]$, whose first projection (on $(0,1]$) is not Borel.
Consider the mapping $\mathbf F$ with domain $[0,1]\times[0,1]$, and signature $\sigma=(f,M)$ (where $f$ is the function symbol and $M$ 
	is a unary relation), with $M(\mathbf F)=S$ and 
$$
f_{\mathbf F}(x,y)=\begin{cases}
	(x,0)&\text{if }y\neq 0\\
	(0,0)&\text{otherwise}
\end{cases}
$$
The mapping $\mathbf M$ is obviously Borel, but fails to be a modeling, as the set $f_{\mathbf F}(S)$ is first-order definable but not Borel.
\end{example}

\begin{definition}
\label{def:pairing}
Let $\mathbf F$ be a Borel $\sigma$-structure with associated probability measure $\nu_{\mathbf F}$, and let
 $\phi\in {\rm FO}(\sigma)$ be a formula with $p$ free variables, such that $\phi(\mathbf F)$ is a Borel subset of $F^p$.
 
 The {\em Stone pairing} of $\phi$ and $\mathbf F$ is the satisfaction probability of $\phi$ in $\mathbf F$ for independent random assignments of elements of $F$ to the free variables of $\phi$ with probability distribution $\nu_{\mathbf F}$, that is:
\begin{equation}
\label{eq:nu}
\langle\phi,\mathbf F\rangle=\nu_{\mathbf F}^{\otimes p}(\phi(\mathbf F)),
\end{equation}
where $\nu_{\mathbf F}^{\otimes p}$ stands for the product measure
$\overbrace{\nu_{\mathbf F}\otimes\dots\otimes\nu_{\mathbf F}}^{p  \text{ times}}$ on $F^p$.
\end{definition}	

Note that if $\mathbf F$ is finite (meaning that $F$ is finite) it holds that
$$
\langle\phi,\mathbf F\rangle=\frac{|\phi(\mathbf F)|}{|F|^p}.
$$

\begin{definition}
Given a fragment $X$ of ${\rm FO}(\sigma)$, a sequence $(\mathbf F_n)_{n\in\bbbn}$ of finite $\sigma$-structures is {\em $X$-convergent} if, for every $\phi\in X$ the limit $\lim_{n\rightarrow\infty}\langle\phi,\mathbf F_n\rangle$ exists.

Moreover, a modeling $\mathbf L$ is a modeling {\em $X$-limit} of the sequence $(\mathbf F_n)_{n\in\bbbn}$ and we note
$\mathbf F_n\xrightarrow{X} \mathbf L$ if, for every first-order formula $\phi\in X$ it holds that
$$
\langle\phi,\mathbf L\rangle=\lim_{n\rightarrow\infty}\langle\phi,\mathbf F_n\rangle.
$$
\end{definition}
Note that if $\mathbf L$ is a modeling $X$-limit of $(\mathbf F_n)_{n\in\bbbn}$, the pairing $\langle\phi,\mathbf L\rangle$ is defined for every first-order formula $\phi$, but its value is required to be equal to  $\lim_{n\rightarrow\infty}\langle\phi,\mathbf F_n\rangle$ only when $\phi$ is in $X$.

Given a fragment $X$ of ${\rm FO}(\sigma)$ (closed under $\vee, \wedge$, and $\neg$)  the equivalence classes of $X$ for logical equivalence form an at most countable Boolean algebra, the Lindenbaum-Tarski algebra $\mathcal L_X$ of $X$. The Stone dual to this algebra is denoted by $S(\mathcal L_X)$. This is a Polish space, the clopen sets of which are in bijection with the elements of $\mathcal L_X$, the topology of which is generated by its clopen sets, and the points of which are the maximal consistent subsets of $\mathcal L_X$ (that is Boolean algebra homomorphisms from $\mathcal L_X$ to the $2$ elements Boolean algebra). 
For instance, if $X={\rm FO}_1^{\rm local}$ then $S(\mathcal L_X)$ is the space of local types.
Considering the Borel $\sigma$-algebra gives $S(\mathcal L_X)$ the structure of a standard Borel space. 

The following representation theorem  was proved in \cite{CMUC}:
\begin{theorem}
\label{thm:mu}
To every finite $\sigma$-structure or $\sigma$-modeling $\mathbf F$ corresponds a unique probability measure $\mu_{\mathbf F}$ on $S(\mathcal L_X)$, such that for every formula $\phi\in X$ it holds that
\begin{equation}
\label{eq:mu}
\langle\phi,\mathbf F\rangle=\int_{S(\mathcal L_X)} I_{K(\phi)}(t)\,{\rm d}\mu_{\mathbf F}(t),
\end{equation}

where $I_{K(\phi)}$ denotes the indicator function of the clopen subset $K(\phi)$ of $S(\mathcal L_X)$ dual to $\phi$. Moreover, a sequence $(\mathbf F_n)_{n\in\bbbn}$ is $X$-convergent if and only if the corresponding sequence of probability measures on $S(\mathcal L_X)$ is weakly convergent. 
\end{theorem}
Note that if the fragment $X$ includes all the fragment ${\rm FO}_0$ of all first-order sentences the support of $\mu_{\mathbf F}$ projects into a single point ${\rm Th}(\mu_{\mathbf F})$ of $S(\mathcal L_{{\rm FO}_0})$, which is (equivalently) characterized by the property
\begin{equation}
	\label{eq:thmu}
	\forall t\in {\rm Supp}(\mu_{\mathbf F})\quad {\rm Th}(\mu_{\mathbf F})=t\cap {\rm QF}_0.
\end{equation}
We call ${\rm Th}(\mu_{\mathbf F})$ the {\em complete theory of $\mu_{\mathbf F}$}, as this is nothing but the complete theory of $\mathbf F$ retrieved from $\mu_{\mathbf F}$.

In this paper we shall be particularly interested by the probability measures  $\mu_{\mathbf F}^{{\rm loc}}$ defined by a $\sigma$-structure $\mathbf F$
on the space $\mathcal T_\infty(\sigma)$ of local types (which is dual to the Lindenbaum-Tarski algebra of local formulas with a single free variable) and
 $\mu_{\mathbf F}^{{\rm loc}(r)}$ defined by a $\sigma$-structure $\mathbf F$
on the (finite) space $\mathcal T_r(\sigma)$ of rank $r$ local types (which is dual to the Lindenbaum-Tarski algebra of local formulas with a single free variable and local rank at most $r$). 

We denote by $\pi_r$ the projection from the space of 
consistent subsets of ${\rm FO}_1^{\rm local}$ to the space of consistent subsets of ${\rm FO}_1^{\rm local}$ with maximum quantifier rank at most $r$.
$$
\pi_r(t)=\{\phi\in t:\ \lrank(\phi)\leq r\}.
$$

Note that $\pi_r$ maps local types to local types with local rank at most $r$.

The mapping $t\mapsto \pi_r(t)$ is measurable and it is immediate that $\mu_{\mathbf F}^{{\rm loc}(r)}$ is the pushforward $\pi_r^*(\mu_{\mathbf F}^{{\rm loc}})$ by $\pi_r$ of the probability measure
$\mu_{\mathbf F}^{{\rm loc}}$ (and that a similar statement holds with any of the probability measures $\mu_{\mathbf F}^{{\rm loc}(r')}$ with $r'>r$).

For an integer $r$ and a $\sigma$-modeling $\mathbf F$, the following easy consequence of \eqref{eq:nu} and \eqref{eq:mu} will be helpfull: for every $t\in \mathcal T_r(\sigma)$ it holds that
\begin{equation}
\label{eq:phit}
\mu_{\mathbf F}^{{\rm loc}(r)}(t)=\nu_{\mathbf F}(\varphi_t(\mathbf F))=\langle\varphi_t,\mathbf F\rangle.
\end{equation}

\subsection{Measuring Proximity}
The topology of ${\rm FO}$-convergence can be metrized by using the following ultrametric
\begin{align}
{\rm d}_{\rm FO}(\mathbf M,\mathbf N)&=\sum_{p\geq 0}\sum_{r\geq 0}
2^{-(p+r)}\dist{p}{r}(\mathbf M,\mathbf N),
\intertext{where}	
\dist{p}{r}(\mathbf M,\mathbf N)&=\sup\Bigl\{|\langle\phi,\mathbf M\rangle-\langle\phi,\mathbf N\rangle|:\ \phi\in{\rm FO}_p, \qrank(\phi)\leq r\Bigr\}.
\end{align}

The following lemma is a direct consequence of \cite[Theorem 13]{limit1}, which in turn follows from Gaifman locality theorem.
\begin{lemma}
\label{lem:locred}
	A mapping modeling $\mathbf L$ is the ${\rm FO}_p$-limit of a sequence of finite mappings if and only if it is both the ${\rm FO}_p^{\rm local}$-limit of a sequence of finite mappings
	and the elementary limit of a sequence of finite mappings.
\end{lemma}

For elementary convergence, the appropriate notion of proximity is the notion of $r$-equivalence, and it holds that $\dist{0}{r}(\mathbf M,\mathbf N)=0$  if and only  $\mathbf M\equiv_r\mathbf N$.

	For local convergence,  we define the following distances (for integers $p\geq 1$ and $r\geq 0$):
\begin{equation}
	\ldist{p}{r}(\mathbf M,\mathbf N)=\sup\Bigl\{|\langle\phi,\mathbf M\rangle-\langle\phi,\mathbf N\rangle|:\ \phi\in{\rm FO}_p^{\rm local}, \lrank(\phi)\leq r\Bigr\}.
\end{equation}	

Note that (by Theorem~\ref{thm:mu}) this is nothing but twice the total variation distance between the probability measures defined by $\mathbf M$ and $\mathbf N$ on the Stone dual of the algebra of local formulas with free variables within $x_1,\dots,x_p$ and local rank at most $r$.

The following lemma is a direct consequence of Lemma~\ref{lem:locred}.
\begin{lemma}
\label{lem:locred2}
	For every fixed signature $\sigma$, every integers $p,r$, and every positive real $\epsilon>0$
	there exist an integer $r'$ and a positive real $\epsilon'>0$, such that for every $\sigma$-modelings $\mathbf M,\mathbf N$ it holds
\begin{equation}
	\mathbf M\equiv_{r'}\mathbf N\text{ and }\ldist{p}{r'}(\mathbf M,\mathbf N)<\epsilon'\quad\Longrightarrow\quad\dist{p}{r}(\mathbf M,\mathbf N)<\epsilon.
\end{equation}
\end{lemma}

In sufficiently sparse structures, where the probability that two random elements are close is small, we can further reduce the computation of the local distance to the case of local formulas with a single free variable:
\begin{lemma}
\label{lem:res}
Let $\delta_{r}(x_1,x_2)$ be the formula ${\rm dist}(x_1,x_2)\leq r$. Then for every integers $p,r$ and every modelings $\mathbf M,\mathbf N$ it holds
\begin{equation}
	\ldist{p}{r}(\mathbf M,\mathbf N)\leq 2p\,\ldist{1}{r}(\mathbf M,\mathbf N)+\binom{p}{2}\bigl(
	\langle\delta_{2r},\mathbf M\rangle+\langle\delta_{2r},\mathbf N\rangle\bigr).
\end{equation}
\end{lemma}
\begin{proof}
Let $\phi$ be a local formula with local rank at most $r$.
The satisfaction of $\phi$ only depends on the $r$-neighborhood of the free variables $x_1,\dots,x_p$. It follows that there exists a finite family $\mathcal F\subseteq\mathcal T_r^p$ such that if ${\rm dist}(v_i,v_j)>2r)$ for every $1\leq i<j\leq p$ then it holds that
$$\mathbf M\models\phi(v_1,\dots,v_p)\quad\iff\quad\mathbf M\models\widehat\phi(x_1,\dots,x_p),$$

where $\widehat\phi$ is the local formula
$$\widehat\phi(x_1,\dots,x_p):=\bigvee_{(t_1,\dots,t_p)\in\mathcal F}\bigwedge_{i=1}^p\varphi_{t_i}(x_i).$$
Moreover, $\widehat\phi(\mathbf M)$ only differs from
$\bigcup_{(t_1,\dots,t_p)\in\mathcal F}\prod_{i=1}^p\varphi_{t_i}(\mathbf M)$ on tuples $(v_1,\dots,v_p)$ with ${\rm dist}(v_i,v_j)\leq 2r$ for some $1\leq i<j\leq p$.
It follows that
$$
\Bigl|\langle\phi,\mathbf M\rangle-\sum_{(t_1,\dots,t_p)\in\mathcal F}\prod_{i=1}^p\langle\phi_{t_i},\mathbf M\rangle\Bigr|<\binom{p}{2}\langle \delta_{2r},\mathbf M\rangle,
$$
as the probability that two random elements of $M$ are at distance at most $2r$ is bounded (by union bound) by $\binom{p}{2}$ times the probability that two random elements are at distance at most $2r$, that is by the right hand side of the inequality.

Of course, the same holds for the modeling $\mathbf N$.

Let $\mu_{\mathbf M}$ (resp. $\mu_{\mathbf N}$) be the probability measure defined by $\mathbf M$ (resp. $\mathbf N$) on $\mathcal T_r(\sigma)$.  
As
$$\sum_{(t_1,\dots,t_p)\in\mathcal F}\prod_{i=1}^p\langle\phi_{t_i},\mathbf M\rangle=\mu_{\mathbf F}^{\otimes p}(\mathcal F),$$
and as it is well known that if $\rho,\lambda$ are probability measures on a finite set it holds that
$$\|\rho^{\otimes p}-\lambda^{\otimes p}\|_{\rm TV}\leq p \|\rho-\lambda\|_{\rm TV}$$
we deduce
\begin{align*}
\frac{1}{2}\Bigl|\sum_{(t_1,\dots,t_p)\in\mathcal F}\prod_{i=1}^p\langle\phi_{t_i},\mathbf M\rangle-\sum_{(t_1,\dots,t_p)\in\mathcal F}\prod_{i=1}^p\langle\phi_{t_i},\mathbf N\rangle\Bigr|
&\leq\|\mu_{\mathbf M}^{\otimes p}-\mu_{\mathbf N}^{\otimes p}\|_{\rm TV}\\
&\leq p\,\ldist{1}{r}(\mathbf M,\mathbf N).
\end{align*}
The statement of the lemma follows.
\end{proof}

\subsection{The Finitary Mass Transport Principle}
\label{sec:fmtp}
The domain of a mapping $\mathbf F$ is partitioned into countably many subsets 
$$F_i=\{x\in F: |f_{\mathbf F}^{-1}(x)|=i\}$$
for $i=0,1,\dots,$ and 
$$F_\infty=\{x\in F: |f_{\mathbf F}^{-1}(x)|=\infty\}.$$

The mass transport principle for mappings takes the following form.

\begin{definition}
The {\em Finitary Mass Transport Principle} (FMTP) for $\mathbf F$
is the satisfaction of the following conditions:
\begin{itemize}
\item $\nu_{\mathbf F}(F_\infty)=0$;
	\item for every measurable subsets $A,B$ of $F\setminus F_\infty$ 
it holds that

\begin{equation}
	\nu_{\mathbf F}(A\cap f_{\mathbf F}^{-1}(B))=
	\int_B|f_{\mathbf F}^{-1}(y)\cap A|\,{\rm d}\nu_{\mathbf F}(y)
\end{equation}
	\end{itemize}
\end{definition}
	
Note that a direct consequence of the FMTP is that for every measurable subset $A$ of $F$ it holds that $\nu_{\mathbf F}(A)\geq \nu_{\mathbf F}(f_{\mathbf F}(A))$.

Intuitively, the FMTP describes the interplay of two measures: the probability measure $\nu_{\mathbf F}$ on $F$ used to randomly select an element, and the counting measure (implicitly) used to count, for instance, the degree of an element. This principle ultimately relies of the fact that the local type of an element is (at least partly) determined by the local type of any of its neighbors. 

\begin{definition}
\label{def:transport}
The {\em transport operator} $\xi$ is a mapping from the space  of consistent subsets of ${\rm FO}_1^{\rm local}$ to itself, defined by
$$
\xi(t)=\{\phi(x)\in{\rm FO}_1^{\rm local}: [(\exists z)\ (z=f(x)\wedge\phi(z))]\in t\}.
$$
\end{definition}

A fundamental property of the transport operator is that if $r'>r$ then for every $\sigma$-structure $\mathbf F$ it holds that
\begin{equation}
\label{eq:xi}
	\Type{\mathbf F}{r}\circ f_{\mathbf F}=\pi_r\circ\xi\circ\Type{\mathbf F}{r'},
\end{equation}
what is depicted by the following diagram:
$$\xy\xymatrix{
F\ar[rr]^{f_{\mathbf F}}\ar[ddd]_{\Type{\mathbf F}{r'}}&&F\ar[d]_{\Type{\mathbf F}{r}}\\
&&\mathcal T_{r}(\sigma)\ar@{.}[r]&\text{rank $r$ local types}\\
&&\xi(\mathcal T_{r'}(\sigma))\ar[u]_{\pi_r}\\
\mathcal T_{r'}(\sigma)\ar[rru]_{\xi}\ar@{.}[rrr]&&&\text{rank $r'$ local types}\\
}\endxy$$


In other words, the rank $r$ local type of the image by $f$ of an element $v$ is exactly the projection of the image by the transport operator of the rank $r+1$ (or any rank $r'>r$) local type of $v$.

We now focus on another aspect of the FMTP.

Let $R>2r$ be positive integers, and
 let $\rho$ be a probability measure on $\mathcal T_R(\sigma)$ (and by extension on $\mathcal T_{r}(\sigma)$). Define
$$T_R(\rho) = \{ \tau \in \mathcal T_{r}:\ \rho(\tau) > 0 \}.$$

For $\tau\in T_{R}(\rho)$ and $t\in T_{r}(\rho)$ define

$$
{\rm adm}^+(\tau,t)=\begin{cases}
	1&\text{if }\varphi_\tau(v)\vdash \varphi_t(f(v))\\
	0&\text{otherwise}
	\end{cases}
$$
and let  ${\rm adm}^-(\tau,t)$
 be the maximum integer $a\in\{0,\dots,{r}+1\}$ such that 
$$\varphi_\tau(v)\vdash \exists x_1,\dots,x_a \Bigl(\bigwedge_{1\leq i\leq a}\bigl(\varphi_t(x_i)\wedge f(x_i)=v\bigr)\ \ \wedge\bigwedge_{1\leq i<j\leq a}(x_i\neq x_j)\Bigr).$$

\begin{definition}
	The probability measure $\rho$ 	satisfies the
	{\em $(R,r)$-restricted FMTP} if there exists a function
	$s:T_{R}(\rho)\times T_{r}(\rho)\rightarrow\{0,1,\dots,r\}\cup ({r},\infty)$, called {\em companion function} of $\rho$, such that for every  $\tau\in T_{R}(\rho)$ and $t\in T_{r}(\rho)$ it holds
\begin{align}
	\min(r,{\rm adm}^-(\tau,t))&=\min(r,s(\tau,t))\\
	\sum_{\tau_1\prec t_1}{\rm adm}^+(\tau_1,t_2)\mu(\tau_1)&=
	\sum_{\tau_2\prec t_2}s(\tau_2,t_1)\mu(\tau_2).
\end{align}
\end{definition}

This notion is justified by the next lemma.

\begin{lemma}
\label{lem:rfmtp}
Let $R>2r$ be positive integers.

Let $\mathbf L$ be mapping modeling $\mathbf L$ that satisfies the FMTP and let $\mu$ be the probability measure on $\mathcal T_R$ defined by $\mu(\tau)=\nu_{\mathbf L}(\varphi_\tau(\mathbf L))$.

Then $\mu$ satisfies the $(R,r)$-restricted FMTP.
\end{lemma}
\begin{proof}
For $\tau\in T_{R}(\mu)$ and $t\in T_{r}(\mu)$ define
\begin{equation}
w(\tau,t)=\frac{\nu_{\mathbf L}(f_{\mathbf L}^{-1}(\varphi_\tau(\mathbf L))\cap \varphi_t(\mathbf L))}{\nu_{\mathbf L}(\varphi_\tau(\mathbf L))}.
\end{equation}

According to FMTP we have the following set of equations
(where $\tau\in T_{R}(\mu)$ and $t\in T_{r}(\mu)$):
\begin{align}
	\min({\rm adm}^-(\tau,t),{r})&=\min(w(\tau,t),{r})\\
	\sum_{\tau_1\prec t_1}{\rm adm}^+(\tau_1,t_2)\mu(\tau_1)&=
	\sum_{\tau_2\prec t_2}w(\tau_2,t_1)\mu(\tau_2).
\end{align}	

\end{proof}
\subsection{The Finite Model Property}

An infinite $\sigma$-structure $\mathbf M$ has the {\em Finite Model Property} if  every sentence $\theta$ satisfied by $\mathbf M$ has a finite model. In other words, $\mathbf M$ has the Finite Model Property if, for every integer $r$, there exists a finite $\sigma$-structure $\mathbf F$ with $\mathbf F\equiv_r \mathbf M$.

Deciding wether an infinite structure has the finite model property is extremely difficult, as deciding wether a sentence has a finite model is undecidable in general, see Trakhtenbrot \cite{Trakhtenbrot1950}.

However, it is clear from our definition that if a modeling $\mathbf M$ is an ${\rm FO}$-limit of a sequence of finite structures then $\mathbf M$ does have the finite model property. When considering the problem of constructing an {\rm FO}-approximation of a modeling $\mathbf M$, we will not only assume that the modeling $\mathbf M$ has the finite model property, but that we can ask an oracle to provide us (for each integer $r$) with a finite structure $\mathbf F$ that is $\mathbf F\equiv_r\mathbf M$.

In some very particular cases, deciding whether a structure has the finite model property and constructing an elementary approximation can be easy. For instance, Lemma~\ref{lem:fmpheight} below asserts that every mapping with finite height has the finite model property and describes how to construct an elementary approximation.
The case of mappings is intermediate between the case of bounded height trees (which have the finite model property) and the case of relational structures with at least one relation symbol with arity at least two, for which satisfiability problem is undecidable.
The {\em Rabin class} $[\text{\em all}, (\omega),(1)]_=$ of first-order logic with equality, one unary function and monadic predicates does not have the finite model property. (For instance, one can consider a sentence expressing that there exists a unique element which is not the image of another element, but that every other element is the image of exactly one element.)
However,  satisfiability problem and finite satisfiability problem for Rabin class are both decidable, though with huge complexity (the first-order theory of one unary function is not elementary recursive). For a general discussion on classical decision problems we refer the reader to \cite{borger2001classical}.
\subsection{Derived Modelings}

Let $\mathbf F$ be a modeling mapping and let $X$ be a non-zero measure first-order definable subset of $F$. We denote by $\restr{\mathbf F}{X}$
the restriction of $\mathbf F$, which is the modeling mapping with domain $X$, probability measure $\nu_{\restr{\mathbf F}{X}}=\frac{1}{\nu_{\mathbf F}(X)}\nu_{\mathbf F}$ and
$$
f_{\restr{\mathbf F}{X}}(v)=\begin{cases}
	f_{\mathbf F}(v)&\text{if }f_{\mathbf F}(v)\in X\\
	v&\text{otherwise}
\end{cases}
$$

\begin{remark}
	The condition that $X$ is first-order definable ensures that $\restr{\mathbf F}{X}$ is a modeling. The condition that $X$ is a Borel subset of $F$ would not be sufficient:
	Consider the modeling mapping $\mathbf F$ with $F=[0,1]\times[0,1]$ and $f_{\mathbf F}$ maps $(x,y)$ to $(x,0)$, and $\nu_{\mathbf F}$ be the usual measure. Then $\mathbf F$ is clearly a modeling.
	 Let $X_0$ be a Borel subset of $(0,1)\times (0,1)$ such that $f_{\mathbf F}(X)$ is not a Borel subset of $[0,1]\times\{0\}$ (such a set can be derived from a standard example of non-Borel $\Sigma_1^1$ sets),
	  and let $X=X_0\cup [0,1]\times\{0\}$. Then $\restr{\mathbf F}{X}$ is not a modeling as the definable subset $\{v: (\exists x)\ (x\neq v)\wedge(f(x)=v)\}$ is not Borel.
\end{remark}

\begin{lemma}
\label{lem:restr}
	Let $\mathbf F$ be a mapping modeling and let $X$ be a non zero-measure first-order definable subset of $F$.
	If $\mathbf F$ satisfies the FMTP then so does $\restr{\mathbf F}{X}$.
\end{lemma}
\begin{proof}
Let $A,B$ be Borel subsets of $X$.
Let $Z=\{v\in X: f_{\mathbf F}(v)\notin X\}$. As $\mathbf F$  satisfies the FTMP it holds
\begin{align*}
	\nu_{\restr{\mathbf F}{X}}(A\cap f_{\restr{\mathbf F}{X}}^{-1}(B))&=
	\nu_{\restr{\mathbf F}{X}}(A\cap f_{\restr{\mathbf F}{X}}^{-1}(B\setminus Z))+
	\nu_{\restr{\mathbf F}{X}}(A\cap f_{\restr{\mathbf F}{X}}^{-1}(B\cap Z))\\
	&=\frac{1}{\nu_{\mathbf F}(X)}
\bigl(\nu_{{\mathbf F}}(A\cap f_{{\mathbf F}}^{-1}(B\setminus Z))+
	\nu_{{\mathbf F}}(A\cap B\cap Z)\bigr)\\
	&=\frac{1}{\nu_{\mathbf F}(X)}
\left(\int_{B\setminus Z}|f_{\mathbf F}^{-1}(y)\cap A|\,{\rm d}\nu_{\mathbf F}(y)+
	\nu_{{\mathbf F}}(A\cap B\cap Z)\right)\\
	&=\int_{B}|f_{\restr{\mathbf F}{X}}^{-1}(y)\cap A|\,{\rm d}\nu_{\restr{\mathbf F}{X}}(y)
\end{align*}
Thus the FMTP holds for $\restr{\mathbf F}{X}$.
\end{proof}

We also note the following:
\begin{lemma}
	\label{lem:mark}
	Let $\mathbf M$ be a modeling and let $\mathbf M^+$ be obtained from $\mathbf M$ by marking exactly one element of $M$ with a new unary relation. Then
	\begin{enumerate}
		\item $\mathbf M^+$ is a modeling;
		\item $\mathbf M^+$ satisfies the FMTP if and only if $\mathbf M$ satisfies the FMTP;
		\item $\mathbf M^+$ has the finite model property if and only if $\mathbf M$ has the finite model property.
	\end{enumerate}
\end{lemma}
\begin{proof}
The first item was proved in \cite{CMUC}.
The second item is obvious as $\mathbf M$ and $\mathbf M^+$ have the same Gaifman graph.
As $\mathbf M$ is a trivial interpretation of $\mathbf M^+$, the finite model property for $\mathbf M^+$ implies the finite model property for $\mathbf M$. Conversely, assume $\mathbf F\equiv_{r+1}\mathbf M$ and start  a Ehrenfeucht-Fra{\"\i}ss{\'e} game by choosing the element that is marked in $\mathbf M^+$. Assume Duplicator follows a winning strategy for the $(r+1)$-rounds game, and mark the vertex chosen by Duplicator in $\mathbf F$. Then (continuing the game) we get that the marked structure is $r$-equivalent to $\mathbf M^+$.
\end{proof}
\subsection{List of Symbols}
Here is a list of the main symbols defined in this section.
\medskip

 \begin{longtable}{|c|p{.8\textwidth}|}
 \hline
Symbol&\multicolumn{1}{c|}{Signification}\\
\hline\hline \endfirsthead
 \hline
Symbol&\multicolumn{1}{c|}{Signification}\\
\hline\multicolumn{2}{c}{}\endhead
\multicolumn{2}{|c|}{Introduced in Section~\ref{sec:prelim_model}}\\
\hline
$\sigma$&signature\\
$\mathbf{F}$&mapping (Definition~\ref{def:mapping})
\\
$F$&domain of structure $\mathbf F$\\
$\phi(\mathbf F)$&set of tuples satisfying $\phi$ in $\mathbf F$\\
$B_r(\mathbf F,u)$&$r$-ball of $u$ in $\mathbf F$\\
\hline
${\rm FO}$&all first-order formulas\\
${\rm FO}_p$&first-order formulas with free variables within $x_1,\dots,x_p$\\
${\rm FO}_0$&sentences\\
${\rm FO}^{\rm local}$&local first-order formulas\\
${\rm FO}_p^{\rm local}$&local first-order formulas with free variables within $x_1,\dots,x_p$\\
${\rm QF}$&quantifier free first-order formulas\\
\hline
${\rm lrank}(\phi)$&local rank of formula $\phi$ (Definition~\ref{def:rank})\\
$t,\tau$&Local types\\
$\mathcal T_r(\sigma)$&set of all rank $r$ local types\\
$\Type{\mathbf F}{r}(v)$&rank $r$ local type of $v$ in $\mathbf F$\\
$\varphi_t(x_1)$&characteristic formula of local type $t$\\
$\delta_r(x_1,x_2)$&formula expressing ${\rm dist}(x_1,x_2)\leq r$\\
$\mathsf I$&interpretation (Definition~\ref{def:interp}) \\
\hline
\hline
\multicolumn{2}{|c|}{Introduced in Section~\ref{sec:lim}}\\
\hline
$\nu_{\mathbf F}$&Probability measure on the domain $F$ of $\mathbf F$ (Definition~\ref{def:modeling}) \\
$\langle\phi,\mathbf F\rangle$&Stone pairing of $\phi$ and $\mathbf F$ (Definition~\ref{def:pairing})\\
$S(\mathcal L_X)$&Stone dual of Lindenbaum-Tarski algebra of $X$\\
$\mu_{\mathbf F}$&Representation measure of  $\mathbf F$ (Theorem~\ref{thm:mu}) \\
$\mu_{\mathbf F}^{\rm loc}$&Representation measure of structure $\mathbf F$ for ${\rm FO}_1^{\rm local}$ fragment\\
${\rm Th}(\mu_{\mathbf F})$&Complete theory of $\mu_{\mathbf F}$\\
$\pi_r$&Projection to consistent subsets of ${\rm FO}_1^{\rm local}$ with quantifier rank at most $r$\\
$\mu_{\mathbf F}^{{\rm loc}(r)}$&Pushforward of $\mu_{\mathbf F}^{\rm loc}$ by $\pi_r$\\
\hline
\hline
\multicolumn{2}{|c|}{Introduced in Section~\ref{sec:fmtp}}\\
\hline
$\zeta$&Transport operator (Definition~\ref{def:transport}) \\
${\rm adm}^+(\tau,t)$&Does $\varphi_\tau(v)$ imply $\varphi_t(f(v))$?\\
${\rm adm}^-(\tau,t)$&How many distinct $u$ with $\varphi_t(u)$ and $f(u)=v$ if  $\varphi_\tau(v)$?\\

\hline
\end{longtable}

\section{First-Order Approximation}
\begin{figure}[ht]
$$\xy\xymatrix@C=5mm@R=5mm{
\text{\small Original modeling}&\mathbf L\ar@{-->}@/_/[d]&&\mathbf F&\text{\small Finite approximation of $\mathbf L$}\\
\text{\small $\epsilon$-residual}&\mathbf L_1\ar@/_/[u]_{\mathsf{I_1}}\ar[dd]\ar@{-->}[dr]&&\mathbf F_1\ar[u]_{\mathsf{I_1}}&\text{\small Finite approximation of $\mathbf L_1$}\\
&&\mathbf E_1\ar`r[ur][ur]&&\text{\small Elementary approximation of $\mathbf L_1$}\\
\text{\small Clean}&\mathbf L_2\ar[d]\ar@{}[r]|{\approx{}_{\Cr{clean}}}&\widetilde{\mathbf L_2}&\mathbf F_2\ar@{-}[uu]&\text{\small Finite approximation of $\mathbf L_2$}\\
\text{\small No short circuits}&\mathbf L_3\ar@{~>}[d]\ar[ur]_{\mathsf{I_2}}&&\mathbf F_3\ar[u]_{\mathsf{I_2}}\ar@{~>}[d]&\text{\small Finite approximation of $\mathbf L_3$}\\
\text{\small Stone measure}&\mu\ar@{-->}[rr]&&\widehat{\mu}\ar@/_/@{..>}[u]&\text{\small Rational Stone measure}\\
}\endxy$$
\captionsetup{singlelinecheck=off,textfont=small}
\caption{Strategy for the proof of Theorem~\ref{thm:invfomap}.
}
\label{fig:stratFO}
\end{figure}
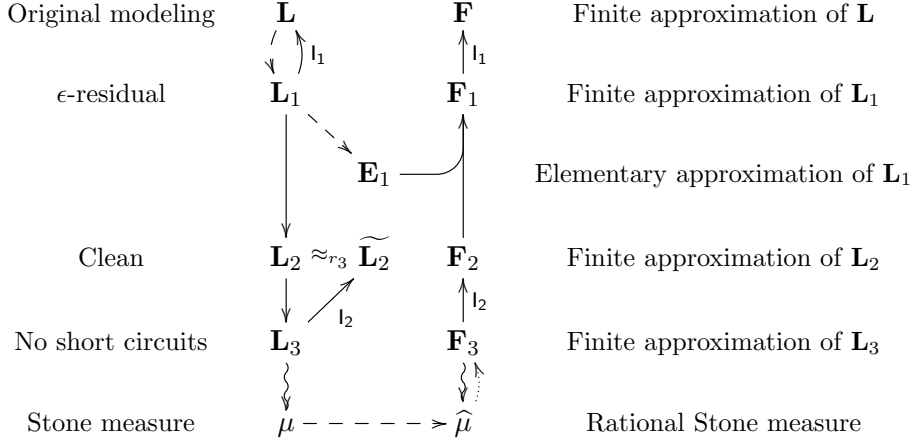

The aim of this section is to prove Theorem~\ref{thm:invfomap}.
The general strategy of the proof is depicted in Fig.~\ref{fig:stratFO}:
\begin{enumerate}
	\item Reduction $\mathbf L\rightarrow\mathbf L_1$, where $\mathbf L_1$  is $\Cr{res}$-residual  (i.e. has no connected component of measure greater than $\Cr{res}$), with recovery interpretation $\mathsf I_1$.
	\item restriction $\mathbf L_1\mapsto\mathbf L_2$ to  no zero-measure rank-$\Cr{clean}$ local types.
	\item Transformation $\mathbf L_2\mapsto\mathbf L_3$  killing all short circuits.
	Interpretation $\mathbf L_3\xrightarrow{\mathsf I_2}\widetilde{\mathbf L_2}$, with local statistics close to $\mathbf L_2$.
	\item Approximation of the rank-$R$ local type measure $\mu$ of $\mathbf L_3$ by a rational measure $\widehat{\mu}$, still satisfying mass transport principle.
	\item Construction of an exact model $\mathbf F_3$ of $\widehat{\mu}$, providing a finite approximation $\mathbf F_3$ of $\mathbf L_3$.
	\item\label{step:rewiring} Rewiring the short cycles by means of interpretation $\mathsf I_2$, leading to an approximation  $\mathbf F_2$ of $\mathbf L_2$.
	\item\label{step:elem} Construction of an elementary approximation $\mathbf E_1$ of $\mathbf L_1$.
	\item\label{step:merge} Merge of $\mathbf E_1$ with a great number of copies of $\mathbf F_2$ to form an ${\rm FO}$-approximation $\mathbf F_1$ of $\mathbf L_1$.
	\item Interpretation $\mathbf F_1\xrightarrow{\mathsf I_1}\mathbf F$ to get an ${\rm FO}$-approximation of the original mapping modeling $\mathbf L$.
\end{enumerate}

We shall reduce the complexity of the approximation problem by requiring more and more properties on the mapping modeling we want to approximate. The different properties we will consider for our mapping modeling are:
\begin{enumerate}[label=(P\arabic*),labelindent=2\parindent, leftmargin=3\parindent, align=left]
	\item\label{it:P1} the modeling measure is atomless;
	\item\label{it:P2} the modeling satisfies the FMTP;
	\item\label{it:P3} the modeling has the finite model property;
	\item\label{it:P4} the modeling is $\Cr{res}$-residual;
	\item\label{it:P5} the modeling is $\Cr{clean}$-clean;
	\item\label{it:P6} the modeling has no cycle of length smaller than $\Cr{cut}$.
\end{enumerate}

During the reduction process, we shall make use of additional unary relations to keep track of the properties of the original mapping. Therefore we shall consider larger and larger signatures:

\begin{itemize}
	\item[$\sigma$] is the signature of both $\mathbf L$ and $\mathbf F$. This signature contains a single unary function symbol $f$ and (possibly) finitely many unary relation symbols.
	\item[$\sigma_1$] is the signature of $\mathbf L_1$, $\mathbf E_1$, and $\mathbf F_1$. It is obtained by adding to $\sigma$ the unary relation symbols $(A_i)_{1\leq i\leq 2\lceil\Cr{res}^{-1}\rceil}$ and $(B_i)_{1\leq i\leq 2\lceil\Cr{res}^{-1}\rceil}$.
	\item[$\sigma_2$] is the signature of $\mathbf L_2$, $\widetilde{\mathbf L_2}$ and $\mathbf F_2$. It is obtained by adding to $\sigma_1$ unary relations $(R_t)_{t\in\mathcal T_{\Cr{clean}}(\sigma_1)}$.
	\item[$\sigma_3$] is the signature of $\mathbf L_3$ and $\mathbf F_3$. It is obtained by adding to $\sigma_2$ unary relations $(U_i)_{1\leq i\leq \Cr{cut}}$ and unary relations $(T_t)_{t\in\mathcal T_{\Cr{clean}}(\sigma_2)}$.
\end{itemize}

We fix integers $p,r$ and a positive real $\epsilon>0$. Our aim is to construct a finite mapping $\mathbf F$ such that 
$\dist{p}{r}(\mathbf L,\mathbf F)<\epsilon$, that is such that
for every first-order formula $\phi$ with at most $p$ free variables and quantifier rank at most $r$, it holds that
$$
	|\langle\phi,\mathbf L\rangle-\langle\phi,\mathbf F\rangle|<\epsilon.
$$
	
We first reduce the problem by separately considering local first-order formulas and sentences.
It follows from Lemma~\ref{lem:locred2}
 that there exist an integer $\Cl[r]{r}$ and a positive real $\Cl[eps]{eps}>0$ such that if $\mathbf L\equiv_{\Cr{r}}\mathbf F$ and $\ldist{p}{\Cr{r}}(\mathbf L,\mathbf F)<\Cr{eps}$
 then it holds $\dist{p}{r}(\mathbf L,\mathbf F)<\epsilon$. We further require $\Cr{eps}<1/16$.

Let
	$\Cl[r]{rr}=4\Cr{r}^2$,
	$\Cl[r]{clean}=2\Cr{rr}+1$,
$\Cl[r]{cut}=\Cr{clean}!$,
	$\Cl[eps]{res}=\Cr{eps}/p^2$,
	$\Cl[eps]{F1}=\Cr{eps}/4p$,
	$\Cl[eps]{epsc}=\Cl[eps]{epsmu}=\Cr{eps}/4p$,
	$\Cl[N]{away}=2\lceil \Cr{res}^{-1}\rceil$,
	$\Cl[r]{elem}=\Cr{r}\Cr{clean}\Cr{away}|\mathcal T_{\Cr{clean}}(\sigma_3)|$.

\subsection{From $\mathbf L$ to $\mathbf L_1$: Reduction to $\epsilon$-residual case}
\label{sec:L1}
For positive real $\epsilon>0$, a modeling $\mathbf M$ is {\em $\epsilon$-residual} if every connected component of $\mathbf M$ has measure at most $\epsilon$.

We consider a signature augmented by $4\lceil \Cr{res}^{-1}\rceil$ marks $A_1,\dots,A_{2\lceil \Cr{res}^{-1}\rceil}$ and
$B_1,\dots,B_{2\lceil \Cr{res}^{-1}\rceil}$, and the basic interpretation $\mathsf I_1$ defined by
$$\eta(x,y):=\left[(f(x)=y))\wedge\neg\bigvee_{i=1}^{2\lceil \Cr{res}^{-1}\rceil} A_i(x)\right]\vee \bigvee_{i=1}^{2\lceil \Cr{res}^{-1}\rceil} (A_i(x)\wedge B_i(y)).$$

We construct a mapping modeling $\mathbf L_1$ from $\mathbf L$ as follows. 

We start by letting $\mathbf L_1$ be a copy of $\mathbf L$, $j=\lceil\Cr{res}^{-1}\rceil+1$, and we modify $\mathbf L_1$ as follows:
We consider the connected component $\mathbf C_i$ ($1\leq i\leq N\leq 1/\epsilon$) of $\mathbf L_1$ with measure $c_i=\nu_{\mathbf L}(C_i)>\Cr{res}$.
If $\mathbf C_i$ contains a non-trivial cycle, we arbitrarily select a vertex $v$ on it, mark $v$ with mark $A_i$, mark $f_{\mathbf L_1}(v)$ by mark $B_i$, and let $f_{\mathbf L_1}(v)=v$.
For $u\in C_i$ let
$$E(u)=\bigcup_{i\geq 1}f_{\mathbf L_1}^{-k}(u).$$

Suppose there exists $v \in C_i$ s.t. $\nu_{\mathbf L_1}(E(v)) > \Cr{res}$. As 
$$\nu_{\mathbf L_1}(E(v)) = \lim_{k \to \infty} \nu_{\mathbf L_1}\Bigl(\bigcup_{1 \leq i \leq k}f_{\mathbf L_1}^{-k}(u)\Bigr),$$
 there exists some $k$ s.t. 
 $$\sum_{u\in f^{-k}(v)}\nu_{\mathbf L_1}(E(u))=\nu_{\mathbf L_1}\Bigl(E(v) \setminus \bigcup_{1 \leq i \leq k} f_{\mathbf L_1}^{-k}(v)\Bigr) \leq \Cr{res}.$$
  Therefore, there is some $u$ s.t. $\nu_{\mathbf L_1}(E(u))> \Cr{res}$ and $\nu_{\mathbf L_1}(E(x)) \leq \Cr{res}$ for all $x \in f_{\mathbf L_1}^{-1}(u)$.

Note that there exist at most $c_i/\Cr{res}$ elements $u\in C_i$ such that $\nu_{\mathbf L_1}(E(u))\geq\Cr{res}$ and $\nu_{\mathbf L_1}(E(x))<\epsilon$ for every $x\in f_{\mathbf L_1}^{-1}(u)$.
For each such element $u$, denoting $W=f_{\mathbf L_1}(u)$, we mark $u$ by a mark $B_j$, every element in $W$ by mark $A_j$, increase $j$ by one, and
redefine $f_{\mathbf L_1}(w)=w$ for every $w\in W$.
As $W$ is first-order definable with a parameter, the structure $\mathbf L_1$ is still a modeling. Doing this, the component $\mathbf C$ gives rise to (possibly uncountably many) small connected components of measure smaller than $\Cr{res}$, and at most one connected component  with measure $\Cr{res}$. 
At the end of the day, we have used up to $2\lceil\Cr{res}^{-1}\rceil$ pairs of marks $A_i$ and $B_i$, $\mathbf L_1$ is $\Cr{res}$-residual, and $\mathbf L=\mathsf I_1(\mathbf L_1)$.

\begin{lemma}
	$\mathbf L_1$ satisfies the properties \ref{it:P1} to  \ref{it:P4} and $\mathbf L = I_1(\mathbf L_1)$.
\end{lemma}
\begin{proof}
As $\nu_{\mathbf L_1}=\nu_{\mathbf L}$, \ref{it:P1} holds for $\mathbf L_1$.
The satisfaction of the FMTP for $\mathbf L$ obviously implies the satisfaction of the FMTP for $\mathbf L_1$ hence \ref{it:P2} holds for $\mathbf L_1$.

The Finite Model Property for $\mathbf L$ implies the one for $\mathbf L_1$ (thus \ref{it:P3} holds): For $r\in\bbbn$, let $\mathbf F$ be a finite mapping such that $\mathbf F\equiv_{r+2\lceil\Cr{eps}^{-1}\rceil} \mathbf L$. Start a Ehrenfeucht-Fra{\"\i}ss{\'e} game of length $r+2\lceil\Cr{eps}^{-1}\rceil$ by selecting in $\mathbf L$ the elements $v_1,\dots,v_N$ marked $B_1,\dots,B_N$ ($N\leq 2\lceil\Cr{eps}^{-1}\rceil$) in $\mathbf L_1$, and let $z_1,\dots,z_N$ be the corresponding elements of $F$ chosen by Duplicator. We construct $\mathbf F_1$ from $\mathbf F$ by marking $z_i$ by mark $B_i$, by marking every element in $Y_i=f_{\mathbf F}^{-1}(z_i)$ by mark $A_i$ and letting $f_{\mathbf F_1}(y)=y$ for every $y\in Y_i$ (for $1\leq i\leq N$). Then it is easily checked that Duplicator's winning strategy for the remaining $r$ steps of the  Ehrenfeucht-Fra{\"\i}ss{\'e} game between $\mathbf L$ and $\mathbf F$ defines a winning strategy for the $r$-step 
 Ehrenfeucht-Fra{\"\i}ss{\'e} game between $\mathbf L_1$ and $\mathbf F_1$ hence $\mathbf F_1\equiv_{r} \mathbf L_1$.
  
 Property \ref{it:P4} holds by construction, as well as the property that $\mathbf L=\mathsf I_1(\mathbf L_1)$.
\end{proof}

\subsection{From $\mathbf L_1$ to $\mathbf L_2$: Cleaning-up}
\begin{definition}
	Let $r\in\bbbn$. A mapping modeling 
	$\mathbf L$ is {\em $r$-clean} if, for every formula $\phi\in{\rm FO}_1^{\rm local}$ with rank at most $r$ it holds that
	$$
\mathbf L\models (\exists x)\phi(x)\quad\iff\quad\langle\phi,\mathbf L\rangle>0.	
$$
\end{definition}

In other words, a mapping modeling 
	$\mathbf L$ is $r$-clean if every local type realized in $\mathbf L$ occurs with non zero probability.

We have proved that $\mathbf L_1$ satisfies \ref{it:P1} to \ref{it:P4}.
We now construct $\mathbf L_2$.

Define
$$T=\{t\in\mathcal T_{\Cr{clean}}(\sigma_1): \langle \varphi_t,\mathbf L_1\rangle>0\},$$
 let $X=\bigvee_{t\in T}\varphi_t(\mathbf L_1)$
--- that is $X$ is the subset of elements of $\mathbf L_1$ whose $\Cr{clean}$-local type appears in $\mathbf L_1$ with no zero probability --- and let $\mathbf L_2$ be obtained from $\restr{\mathbf L_1}{X}$ by  adding marks $R_t$ ($t\in \mathcal T_{\Cr{clean}}(\sigma_1)$), in such a way that for all $t\in \mathcal T_{\Cr{clean}}(\sigma_1)$ and $v\in L_2$ it holds that
$$
\mathbf L_2\models R_t(v)\quad\iff\quad \mathbf L_1\models\varphi_t(v)\quad\iff\quad\Type{\mathbf L_1}{\Cr{clean}}(v)=t.
$$

\begin{lemma}
	The mapping modeling $\mathbf L_2$ satisfies properties \ref{it:P1} to \ref{it:P5}.
\end{lemma}
\begin{proof}
Let $\widehat{\mathbf L}_1$ be the $\sigma_2$-mapping obtained by the trivial interpretation adding marks $R_t$ in such a way that
$R_t(\widehat{\mathbf L}_1)=\varphi_t(\mathbf L_1)$.
As we made use of a trivial interpretation, $\widehat{\mathbf L}_1$ is a modeling and properties \ref{it:P1} to \ref{it:P4} still hold. Note that
$\mathbf L_2=\restr{\widehat{\mathbf L}_1}{X}$.
It is immediate that \ref{it:P1} and \ref{it:P4} hold.
According to Lemma~\ref{lem:restr}, \ref{it:P2} holds.
If $\mathbf F$ is a finite elementary approximation of $\widehat{\mathbf L}_1$ then $\restr{\mathbf F}{X}$ 
 is a finite elementary approximation of $\mathbf L_2$ hence $\mathbf L_2$ has the finite model property \ref{it:P3}.
An easy $\Cr{clean}$-step local Ehrenfeucht-Fra{\"\i}ss{\'e} game easily shows that if $u,v\in L_2$ have same rank $\Cr{res}$ local type in $\mathbf L_1$ then they have the same rank $\Cr{clean}$ local type in $\mathbf L_2$. It follows that $\mathbf L_2$ is $\Cr{clean}$-clean thus \ref{it:P5} holds.
\end{proof}

\subsection{From $\mathbf L_2$ to $\mathbf L_3$: Cutting the short cycles}
\begin{figure}
	\begin{center}
		\includegraphics[width=\textwidth]{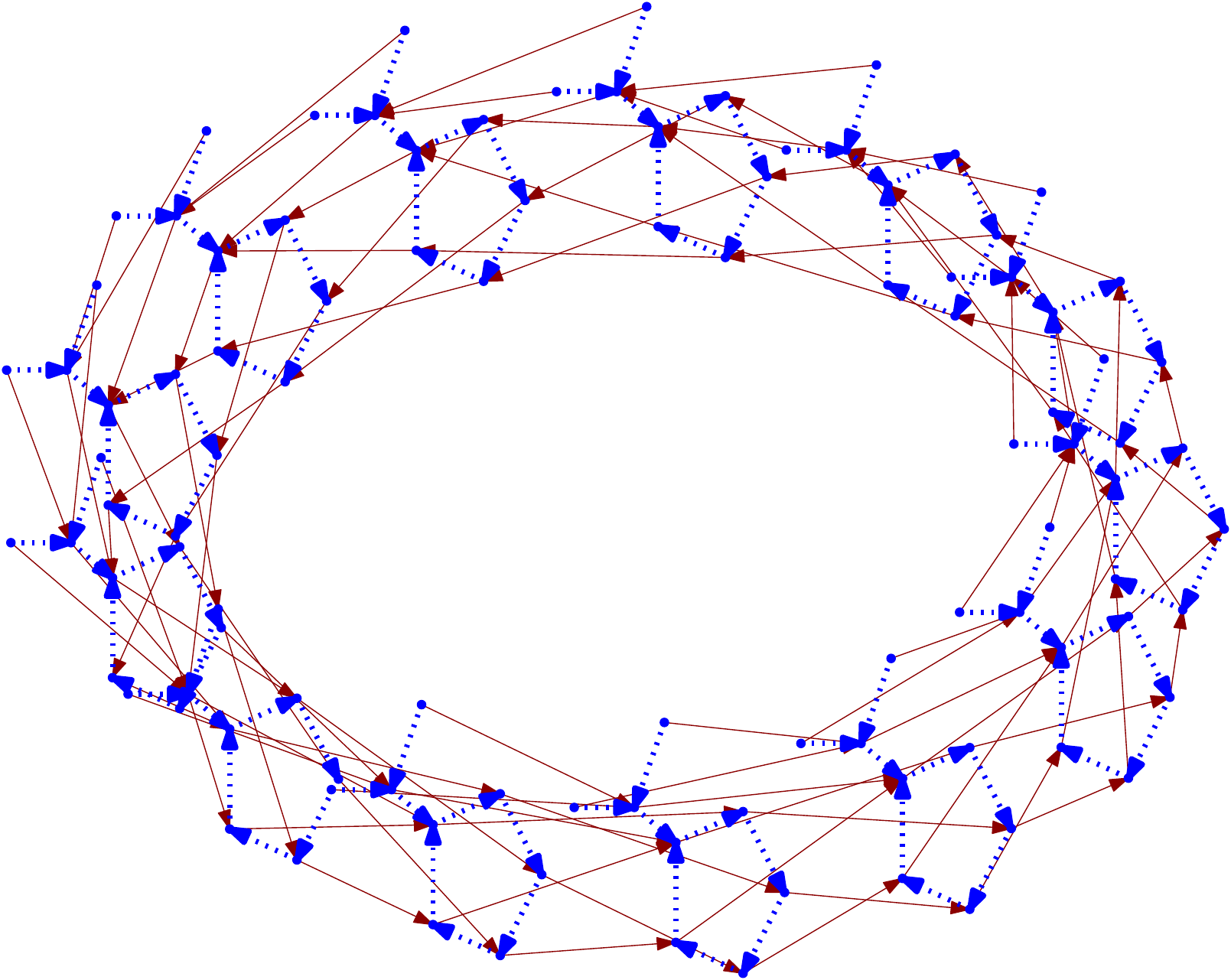}
	\end{center}
	\caption{Construction of ${\mathbf L_3}$}
	\label{fig:killcyc}
\end{figure}

Cutting the short cycles will allow to handle mapping modelings that are locally acyclic, which will strongly simplify the proofs.
A natural procedure would be to consider a Borel transversal of all short cycles (which exists thanks to Borel selection theorem \cite[p. 78]{Kechris1995}), to mark it, and to use an interpretation to kill the cycles at the mark. However, such an approach fails as marking a Borel subset of a modeling does not in general keep the property of being a modeling (see Example~\ref{ex:souslin}).
We shall use a different approach. Let $\Gamma$ be the set $[\Cr{cut}]$. We consider the $\sigma_3$-mapping modeling $\mathbf L_3$ with domain $L_3=L_2\times\Gamma$, measure $\nu_{\mathbf L_3}=\nu_{\mathbf L_2}\otimes\delta_\Gamma$ (where $\delta_\Gamma$ is the uniform measure on $\Gamma$), with $(x,i)$ marked by $U_i$,  $T_{\Type{\mathbf L_2}{\Cr{clean}}(x)}$, 
	and
$$
f_{\mathbf L_3}(x,i)=(f_{\mathbf L_2}(x), i+1 \bmod \Cr{cut}).
$$
An example of construction of $\mathbf L_3$ is shown on Fig.~\ref{fig:killcyc}.

\begin{lemma}
	The mapping modeling $\mathbf L_3$ satisfies \ref{it:P1} to \ref{it:P6}.
\end{lemma}
\begin{proof}
Property \ref{it:P1} obviously holds.

As $\mathbf L_2$ satisfies the FMTP, so does $\mathbf L_3$. Indeed, let $A,B$ be Borel subsets of $L_3$ such that ${\rm deg}_B^{\mathbf L_3}(v)$ is bounded for $v\in A$ and ${\rm deg}_A^{\mathbf L_3}(v)$ is bounded for $v\in B$. Then we can write $A=\bigcup_i A_i\times\{i\}$ and $B=\bigcup_j B_j\times\{j\}$, where the $A_i$'s and the $B_j$'s are Borel subsets of $L_2$. Then it holds that
\begin{align*}
	\nu_{\mathbf L_3}(A\cap f_{\mathbf L_3}^{-1}(B))&=
	\frac{1}{\Cr{cut}}\sum_i \nu_{{\mathbf L}_2}(A_i\cap f_{{\mathbf L}_2}^{-1}(B_{i+1 \bmod \Cr{cut}}))\\
	&=\frac{1}{\Cr{cut}}\sum_j\int_{B_j}|f_{\mathbf L_2}^{-1}(y)\cap A_{j-1\bmod \Cr{cut}}|\, {\rm d}\nu_{\mathbf L_2}(y)\\
	&=\int_{B}|f_{\mathbf L_3}^{-1}(y)\cap A|\, {\rm d}\nu_{\mathbf L_3}(y)
\end{align*}
Hence \ref{it:P2} holds.

It is immediate that if for some $R\in\bbbn$ it holds that $\mathbf F\equiv_R \mathbf L_2$ then
if $\mathbf F'$ is obtained from $\mathbf F$ in the same way that
$\mathbf L_3$ is obtained from $\mathbf L_2$ it holds that $\mathbf F'\equiv_R \mathbf L_3$ (Duplicator's strategy immediately follows from its strategy in an Ehrenfeucht-Fra{\"\i}ss{\'e} game between $\mathbf F$ and $\mathbf L_2$). Thus \ref{it:P3} holds.

It is easily checked that the measure of a connected component of
$\mathbf L_3$ is at most the measure of its projection on $\mathbf L_2$. Thus \ref{it:P4} holds.

As $\Cr{cut}>\Cr{clean}$, an easy Ehrenfeucht-Fra{\"\i}ss{\'e} game shows that if two elements $x,y$ of $L_2$ have the same $\Cr{clean}$ local type in $\mathbf L_2$ and $1\leq i,j\leq $ the $(x,i)$ and $(y,j)$ have the same $\Cr{clean}$ local type in $\mathbf L_3$. Thus, as $\mathbf L_2$ is $\Cr{clean}$-clean so is $\mathbf L_3$.
 Hence \ref{it:P5} holds for $\mathbf L_3$.
 
 By construction, $\mathbf L_3$ has no cycle of length smaller than $\Cr{cut}$ thus \ref{it:P6} holds.
\end{proof}

For $1\leq\ell\leq\Cr{clean}$ let ${\rm Z}_\ell$ be the subset of all the $t\in \mathcal T_{\Cr{clean}}(\sigma_2)$ that contain the formula
$\bigl[(f^\ell(x)=x)\wedge\bigwedge_{i<\ell}(f^i(x)\neq x)\bigr]$ (which means that $x$ belongs to a cycle of length $\ell$).

Now we consider the basic interpretation $\mathsf I_2$, with
$$
\eta(x,y):=\biggl[\bigvee_{\ell=1}^{\Cr{clean}}(\zeta_\ell(x)\wedge(y=f^{\ell-1}(x))\biggr]\vee \biggl[(y=f(x))\wedge \neg\bigvee_{\ell=1}^{\Cr{clean}}U_i(x)\biggr],
$$

where 
$$
\zeta_l(x):=U_l(x)\wedge \bigvee_{t\in {\rm Z}_\ell}T_t(x),
$$
which also forgets the marks $U_i$ and $T_t$.
Let $\widetilde{\mathbf L}_2=\mathsf I_2(\mathbf L_3)$.

\begin{lemma}
	For every $\phi\in{\rm FO}_1^{\rm local}$ with rank at most $\Cr{clean}$ it holds that
$$
\langle\phi,\widetilde{\mathbf L}_2\rangle=\langle\phi,\mathbf L_2\rangle.
$$
\end{lemma}
\begin{proof}
It is straightforward that for every $v\in L_2$ and every $i\in\Gamma$ it holds that
$$
\Type{\widetilde{\mathbf L}_2}{\Cr{clean}}(v,i)=\Type{\mathbf L_2}{\Cr{clean}}(v).
$$
Hence for every $\phi\in{\rm FO}_1^{\rm local}$ with rank at most $\Cr{clean}$ it holds that
$$
\langle\phi,\widetilde{\mathbf L}_2\rangle=\langle\phi,\mathbf L_2\rangle.
$$
\end{proof}

\subsection{From $\mu$ to $\widehat{\mu}$: Approximating the Stone measure}

~\\
Let $\mu=\mu_{\mathbf L_3}^{{\rm loc}(\Cr{clean})}$. As $\mathbf L_3$ satisfies the FMTP, according to Lemma~\ref{lem:rfmtp}, the probability measure $\mu$ satisfies the $({\Cr{clean}},{\Cr{rr}})$-restricted FMTP.
\begin{lemma}
\label{lem:ratmu}
There exists a rational probability measure $\widehat\mu$ on $\mathcal T_{\Cr{clean}}(\sigma_3)$ with same support as $\mu$, that satisfies the $({\Cr{clean}},{\Cr{rr}})$-restricted MTP, and such that $\|\mu-\widehat\mu\|_{\rm TV}<\Cr{epsmu}$.
\end{lemma}
\begin{proof}
Let $w$ be a companion function for $\mu$, and let
\begin{align*}
Q_1&=\{(\tau,t)\in T_{\Cr{clean}}(\mu)\times T_{\Cr{rr}}(\mu): w(\tau,t)\leq {\Cr{rr}}\}\\
Q_2&=\{(\tau,t)\in T_{\Cr{clean}}(\mu)\times T_{\Cr{rr}}(\mu): w(\tau,t)>{\Cr{rr}}\}\\	
\end{align*}
Consider the following set of Diophantine equations and inequalities with variables $x_\tau$ ($\tau\in T_{\Cr{clean}}(\mu)$) and $y_{\tau,t}$ ($(\tau,t)\in Q_2$):

$$x_\tau>0,\quad\sum_{\tau\in T_{\Cr{clean}}}x_\tau= 1,\quad y_{\tau,t}\geq 0,$$
$$\sum_{\tau_1\prec t_1}{\rm adm}^+(\tau_1,t_2)x_{\tau_1}=
	\sum_{\substack{\tau_2\prec t_2\\(\tau_1,t_2)\in Q_1}}{\rm adm}^-(\tau_1,t_2)x_{\tau_2}+\sum_{\substack{\tau_2\prec t_2\\(\tau_1,t_2)\in Q_2}}({\Cr{rr}}x_{\tau_2}+y_{\tau_1,t_2})
	$$
Then this set defines a convex polytope containing a solution for  $x_\tau=\mu(\tau)$ and 
$y_{\tau,t}=(w(\tau,t)-{\Cr{rr}})\mu(\tau)$.

 Since this polytope has rational vertices, either the aforementioned solution is rational, or there is a strictly positive rational solution in any of its neighborhood.
Let $(\widehat x_\tau, \widehat y_{\tau,t})$ be such a rational solution, such that $\sum_{\tau\in T_{\Cr{clean}}(\sigma_3)}|x_\tau-\widehat{x}_\tau|<\Cr{epsmu}$.

Define $\widehat{\mu}(\tau)=\widehat x_\tau$. Then  $\widehat{\mu}$ has same support as $\mu$ and 
$\|\mu-\widehat\mu\|_{\rm TV}<\Cr{epsmu}$,
and $\widehat\mu$, with companion function 
$s(\tau,t)={\Cr{rr}}+\widehat y_{\tau,t}/\widehat x_\tau$, satisfies the $({\Cr{clean}},{\Cr{rr}})$-restricted FMTP. 
\end{proof}
\subsection{Constructing $\mathbf F_3$}

It is possible, by means of a (relatively low local rank) local formula, to specify that in the neighborhood of an element $v$, related in a given way (by means of a digraph $D$ indicating which element is the image of which element), one finds an element $u_1$ with rank $\rho_1$ local type $t_1$, an element $u_2$ with rank $\rho_2$ local type $t_2$,\dots, and  an element $u_k$ with rank $\rho_k$ local type $t_k$. This is the aim of the following definition.

\begin{definition}
Let $\sigma$ be a mapping signature, 
let $k\in\bbbn$, $\rho_1>\rho_2>\dots>\rho_k\geq 0$, $t_1\in\mathcal T_{\rho_1}(\sigma),\dots, t_k\in \mathcal T_{\rho_k}(\sigma)$, and let $D\subseteq [k+1]\times[k+1]$ be the arc set of a digraph with outdegrees at most $1$ and  connected underlying graph.
We define the {\em characteristic formula} $\theta\in {\rm FO}_1^{\rm local}(\sigma)$ of $((\rho_i)_{i\in[k]},(t_i)_{i\in[k]}, D)$
inductively as follows:
\begin{align*}
\theta_{k+1}(x_1, \dots, x_{k+1}) &:= \bigwedge_{1 \leq i < j \leq k+1}(x_i \neq x_j) \wedge \bigwedge_{(i, j) \in D} f(x_i) = x_j\\
\theta_i(x_1, \dots, x_i) &:= \exists y_i\,[\varphi_{t_i}^{\rho_i}(y_i) \wedge \theta_{i+1}(x_1, \dots, x_i, y_i)]&\text{$(1\leq i\leq k)$}\\
\theta(x)&:=\theta_1(x)
\end{align*}
\end{definition}
Note that the rank of $\theta$ is at most $\rho_1+1=\max\{\rho_i+i:\ 1\leq i\leq k\}$.

\begin{lemma}
\label{lem:lef}
	Let $\mathbf F$ be a $\sigma_3$-mapping with no cycle of length $1<\ell\leq \Cr{cut}$, and let $\Upsilon:F\rightarrow\mathcal T_{\Cr{clean}}(\sigma)$ be such that
	\begin{enumerate}
		\item for every unary mark $M$ in the signature and every $v\in F$, 
		  $M(v)$ holds in $\mathbf F$ if and only if $M(x)\in\Upsilon(v)$;
	     \item\label{it:acyc} for every $1\leq\ell\leq \Cr{cut}$ and every $v\in F$ it holds that
	     $[f^i(x)=x]\notin \Upsilon(v)$.
		\item\label{it:propagate} for every $v\in F$ it holds that
	     $${\rm adm}^+(\Upsilon(v),\pi_{\Cr{rr}}(\Upsilon(f_{\mathbf F}(v))))=1;$$
		\item\label{it:back} for every $v\in F$ and $t\in \mathcal T_{\Cr{rr}}(\sigma_3)$ it holds that
	     $$\min\bigl(\Cr{rr},{\rm adm}^-(\Upsilon(v),t)\bigr)=
	     \min\bigl(\Cr{rr},|\{u\in f_{\mathbf F}^{-1}(v): \pi_{\Cr{rr}}(\Upsilon(u))=t\}|\bigr).$$
	\end{enumerate}

Then for every $v\in\mathbf F$ it holds that $\Type{\mathbf F}{\Cr{rr}}(v)=\pi_{\Cr{rr}}(\Upsilon(v))$.
\end{lemma}
\begin{proof}
First note that Property \ref{it:propagate} implies that for every 
$0\leq i\leq \Cr{rr}$ and every $v\in F$ it holds that
$$
(\pi_i\circ \Upsilon)\circ f_{\mathbf F}=\xi\circ(\pi_{i+1}\circ\Upsilon)x.
$$

Note that this is analog to \eqref{eq:xi}, which states that for every non-negative integer $i$ and every mapping $\mathbf M$ it holds that
$$
\Type{\mathbf M}{i}\circ f_{\mathbf M}=\xi\circ \Type{\mathbf M}{i+1}.
$$

For $v\in F$, let $\mathbf M$ be a countable model of 
$(\exists x)\,\varphi_{\Upsilon(v)}(x)$, and let $z\in M$ be such that $\mathbf M_0\models\varphi_{\Upsilon(v)}(z)$, that is $\Type{\mathbf M_0}{\Cr{clean}}(z)=\Upsilon(v)$. By Property \ref{it:propagate} it holds that $f^{d}(x)=x$ belongs to no $\Upsilon(u)$ at distance at most $\Cr{clean}-d$ from $z$.
Considering the ball of radius $\Cr{clean}+1$ around $z$ we deduce that there exists a connected mapping $\mathbf M$ with a special element $z$, which has no cycle of length $>1$ (hence the Gaifman graph of $\mathbf M$ is a tree), at most one fixed point at distance $\Cr{clean}+1$ from $z$, and such that $\Type{\mathbf M}{\Cr{clean}}(z)=\Upsilon(v)$.

In order to prove $\Type{\mathbf F}{\Cr{rr}}(v)=\pi_{\Cr{rr}}\circ\Upsilon(v)=\Type{\mathbf M}{\Cr{rr}}(z)$ it is sufficient to prove that Duplicator has a winning strategy for the $\Cr{rr}$ steps local Ehrenfeucht-Fra{\"\i}ss{\'e} game between $(\mathbf F,v)$ and $(\mathbf M,z)$.

		Assume that for some $0\leq k<\Cr{rr}$ we have
		$v_0,\dots,v_k\in F$ and $z_0,\dots,z_k\in\mathbf M$ with  $v_0=v$ and $z_0=z$, such that $v_i\mapsto z_i$ is a partial isomorphism, and such that
		for every $0\leq i\leq k$ it holds that
		$$
		\Type{\mathbf M}{\Cr{rr}-i}(z_i)=\pi_{\Cr{rr}-i}\circ\Upsilon(v_i).
		$$
		
	Now consider a Spoiler move. There are six cases:
\begin{enumerate}[i)]
	\item\label{case:1} Spoiler chooses $v_{k+1}\in F$, and there exists $0\leq a<k+1$ such that $f_{\mathbf F}(v_a)=v_{k+1}$.
	
	In this case, $\Type{\mathbf M}{\Cr{rr}-a}(z_a)=\pi_{\Cr{rr}-a}\circ\Upsilon(v_a)$  implies
\begin{align*}
	\Type{\mathbf M}{\Cr{rr}-a-1}\circ f_{\mathbf M}(z_a)&=\xi\circ\Type{\mathbf M}{\Cr{rr}-a}(z_a)\\
	&=\xi\circ \pi_{\Cr{rr}-a}\circ\Upsilon(v_a)\\
	&=\pi_{\Cr{rr}-a-1}\circ\Upsilon\circ f_{\mathbf F}(v_a)
\end{align*}
	 Thus we can let $z_{k+1}=f_{\mathbf M}(z_a)$.
	\item\label{case:2} Spoiler chooses $v_{k+1}\in F$, there exists $0\leq a<k+1$ such that $f_{\mathbf F}(v_{k+1})=v_a$, and  for every $0\leq i<a$ it holds $f_{\mathbf F}(v_i)\neq v_a$.
	
	Let $b_1<b_2<\dots<b_{\ell+1}=k+1$ be such that
	$f_{\mathbf F}^{-1}(v_a) \cap\{v_0,\dots,v_{k+1}\}=\{v_{b_1},\dots,v_{b_{\ell+1}}\}$. Note that $b_1>a$ by assumption.
	
	For $1\leq i\leq \ell+1$, let $\rho_i=\Cr{rr}-b_i$ and $t_i=\pi_{\rho_i}\circ\Upsilon(v_{b_i})$. Let $D$ be the set of pairs $(i,1)$ for $2\leq i\leq \ell+2$, and let $\theta(x)$ be the characteristic formula of $((\rho_i)_{i\in [\ell+1]}, (t_i)_{i\in [\ell+1]}, D)$. This formula has rank at most $\rho_1+1\leq \Cr{rr}-a$ so it holds that
$\theta(x)\in\pi_{\Cr{rr}-a}\circ\Upsilon(v_a)=\Type{\mathbf M}{\Cr{rr}-a}(z_a)$.
 Thus there exist
$z_{b_1}',\dots,z_{b_\ell}',z_{k+1}'$ in $f_{\mathbf M}^{-1}(z_a)\setminus\{z_a\}$, such that 
\begin{equation}
\label{eq:zprime}
\Type{\mathbf M}{\rho_i}(z_{b_i}')=\pi_{\Cr{rr}-b_i}\circ\Upsilon(v_{b_i})\qquad(1\leq i\leq \ell+1).
\end{equation}
If $z_{b_i}'$ is not equal to $z_{b_i}$ for every $1\leq i\leq \ell$, let $i$ be minimum such that $z_{b_i}'\neq z_{b_i}$.

\begin{itemize}
	\item If $z_{b_i}=z_{b_j}'$ for some $j>i$ then it holds that
$$t_j=\Type{\mathbf M}{\rho_j}(z_{b_j'})=\Type{\mathbf M}{\rho_j}(z_{b_i})\subseteq \Type{\mathbf M}{\rho_i}(z_{b_i})=t_i$$
and we deduce that \eqref{eq:zprime} still holds after exchange
of $z_{b_i}'$ and $z_{b_j}'$.
\item Otherwise, we let $z_{b_i}'=z_{b_i}$ and remark that \eqref{eq:zprime} still holds.
\end{itemize}

We repeat this process until we get $z_{b_i}'=z_{b_i}$ for every $1\leq i\leq \ell$.  Then we let $z_{k+1}=z_{k+1}'$. 

	\item\label{case:3} Spoiler chooses $v_{k+1}\in F$, there exists $0\leq a<k+1$ such that $f_{\mathbf F}(v_{k+1})=v_a$, and there exists $0\leq i<a$ such that it holds that $f_{\mathbf F}(v_i)=v_a$.

	Let $b_1 < b_2 < \dots < b_{\ell+1} = k+1$ be such that $f_{\mathbf F}^{-1}(v_a)\cap \{ v_0, \dots, v_{k+1} \} = \{ v_{b_1}, \dots, v_{b_{\ell+1}} \}$.
	
	Note that there can be only one $0 \leq i < a$ s.t. $f_{\mathbf F}(v_i) = v_a$, as otherwise the two vertices would not be connected before step $a$, so $b_1 < a < b_2$. 
		
	For $1 \leq i \leq \ell + 1$, let $\rho_1 = \Cr{rr} - a$, $t_1 = \pi_{\Cr{rr}-a} \circ \Upsilon(v_a)$, and $\rho_i = \Cr{rr} - b_i, t_i = \pi_{\rho_i} \circ \Upsilon(v_{b_i})$ for $2 \leq i \leq \ell + 1$.
	
	Let $D$ be the set of pairs $(i, 2)$ for $i \in \{1, \dots, \ell + 2 \} \setminus \{2\}$, and let $\theta(x)$ be the characteristic formula of $((\rho_i)_{i \in [\ell+1]}, (t_i)_{i \in \ell + 1}, D)$. This formula has rank at most $\rho_1 + 1 \leq \Cr{rr}-a+1 \leq \Cr{rr} - b_1$ so it holds that $\theta(x) \in \pi_{\Cr{rr}-b_1} \circ \Upsilon(v_{b_1}) = \Type{\mathbf M}{\Cr{rr}-b_1}(z_{b_1})$. Thus, there exists $z'_a, z'_{b_2}, \dots, z'_{b_\ell}, z'_{k+1} \in f_{\mathbf M}^{-1} \circ f_{\mathbf M}(z_{b_1})$, all distinct, such that $f_{\mathbf M}(z_{b_1}) = z_a$ and
	
\begin{equation}
\label{eq:zprime}
\Type{\mathbf M}{\rho_i}(z_{b_i}')=\pi_{\Cr{rr}-b_i}\circ\Upsilon(v_{b_i})\qquad(2\leq i\leq \ell+1).
\end{equation}

As in the previous case, we can assume $z'_{b_i} = z_{b_i}$ for $2 \leq i \leq \ell$ and let $z_{k+1} = z'_{k+1}$.

	\item\label{case:4} Spoiler chooses $z_{k+1}\in M$, and there exists $0\leq a<k+1$ such that $f_{\mathbf M}(z_a)=z_{k+1}$.
	
	As in Case~\ref{case:1}, $\Type{\mathbf M}{\Cr{rr}-a}(z_a)=\pi_{\Cr{rr}-a}\circ\Upsilon(v_a)$  implies
	$$\Type{\mathbf M}{\Cr{rr}-a-1}\circ f_{\mathbf M}(z_a)=\pi_{\Cr{rr}-a-1}\circ\Upsilon\circ f_{\mathbf F}(v_a).$$
	 Thus we can let $v_{k+1}=f_{\mathbf F}(v_a)$.
	\item\label{case:5} Spoiler chooses $z_{k+1}\in M$, there exists $0\leq a<k+1$ such that $f_{\mathbf M}(z_{k+1})=z_a$, and  for every $0\leq i<a$ it holds $f_{\mathbf M}(z_i)\neq z_a$.
		
Let  $\tau=\Upsilon(v_a)$, let $t=\Type{\mathbf M}{\Cr{rr}-(k+1)}(z_{k+1})$, and let $p$ be the number of elements of $f_{\mathbf M}^{-1}(z_a)\cap\{z_0,\dots,z_{k+1}\}$ with rank $(\Cr{rr}-(k+1))$ local type $t$.

By assumption, it holds that
$\Type{\mathbf M}{\Cr{rr}-a}(z_a)=\pi_{\Cr{rr}-a}(\tau)$. Thus
$$\sum_{t'\prec t}{\rm adm}^-(\tau,t')\geq p,$$
where the sum is over local types $t'\in \mathcal T_{\Cr{rr}}(\sigma_3)$ such that $t'\prec t$.
According to Property \ref{it:back}, it holds that
$$\sum_{t'\prec t}{\rm adm}^-(\tau,t')=|\{u\in f_{\mathbf F}^{-1}(v_a):\ \pi_{\Cr{rr}-(k+1)}(\Upsilon(u))=t\}|.$$
 It follows that there exists $v_{k+1}\in f_{\mathbf F}^{-1}(v_a)$, distinct from $v_0,\dots,v_k$, such that $\pi_{\Cr{rr}-(k+1)}(\Upsilon(v_{k+1}))=\Type{\mathbf M}{\Cr{rr}-(k+1)}(z_{k+1})$.

	\item\label{case:6} Spoiler chooses $z_{k+1}\in M$, there exists $0\leq a<k+1$ such that $f_{\mathbf M}(z_{k+1})=z_a$, and there exists $0\leq i<a$ such that it holds that $f_{\mathbf M}(z_i)=z_a$.
		
This case is solved similarly, by considering the element $z_i$ such that $f_{\mathbf M}(z_i)=z_a$, and showing that the number of elements of $f_{\mathbf F}^{-1}(v_i)$ with same rank $(\Cr{rr}-(k+1))$ local type as $z_{k+1}$ is at least equal to the number of elements of $f_{\mathbf F}^{-1}(z_i)\cap\{z_0,\dots,z_{k+1}\}$ with same rank $(\Cr{rr}-(k+1))$ local type as $z_{k+1}$.
\end{enumerate}
\end{proof}

\begin{lemma}
Let $\Cr{clean}>2\Cr{rr}$ be positive integers, and let $\widehat\mu$ be a rational probability measure on $\mathcal T_{\Cr{clean}}(\sigma_3)$, such that
\begin{enumerate}
	\item $\widehat\mu$ is clean: for every $\tau\in\mathcal T_{\Cr{clean}}(\sigma_3)$ with $\widehat\mu(\tau)>0$ and for every $t\in\mathcal T_{\Cr{clean}-1}(\sigma_3)$, if
	$\phi_t(f(x))\in \tau$ 
	then $\sum_{\tau'\prec t}\widehat\mu(\tau')>0$;
	\item for every $1<i\leq\Cr{clean}$ the formula
$f^i(x)=x$ does not belong to any $\tau\in \mathcal T_{\Cr{clean}}(\sigma_3)$ with positive $\widehat\mu$-measure;
	\item the measure $\widehat\mu$  satisfies the $(\Cr{clean},\Cr{rr})$-restricted MTP.
\end{enumerate}

Then there exists a finite $\sigma_3$-mapping $\mathbf F_3$ such that 
for every local formula $\phi\in{\rm FO}_1^{\rm local}(\sigma_3)$ with local rank at most $\Cr{rr}$ it holds that
\begin{equation}
\label{eq:F3}	
\langle\phi,\mathbf F_3\rangle=\sum_{\tau\ni\phi}\widehat\mu(\tau).
\end{equation}
\end{lemma}
\begin{proof}
	Let $T_{\Cr{clean}}=\{\tau\in\mathcal T_{\Cr{clean}}(\sigma_3): \widehat\mu(\tau)>0\}$ and $T_{\Cr{rr}}=\{\pi_{\Cr{rr}}(\tau): \tau\in T_{\Cr{clean}}\}$.

Let $\Cl[N]{N}\in\bbbn$ be such that $\Cr{N}\,\widehat\mu$ is integral, 
and let $\zeta:[\Cr{N}]\rightarrow T_{\Cr{clean}}$ be such that for every $\tau\in T_{\Cr{clean}}$ it holds $|\zeta^{-1}(\tau)|=\Cr{N}\,\widehat\mu(\tau)$.

We construct a (partial) mapping $g:[\Cr{N}]\rightarrow [\Cr{N}]$ inductively. 
We start with an empty domain. For each $i\in[N]$ (not yet in the domain), let $t=\pi_r(\zeta(i))$.
We consider the elements of $[\Cr{N}]$ such that ${\rm adm}^-(\zeta(j),t)$
is either $\Cr{rr}+1$, or greater than the number of $k\in g^{-1}(j)$ such that $\zeta(k)\prec t$. Among these elements, we choose one element $j$ such that ${\rm adm}^-(t,\zeta(j)$ is minimal, and let $g(i)=j$.

Now we prove that the above construction never gets stuck and that, at the end of the day, for every $j\in[\Cr{N}]$ and every $t\in\mathcal T_{\Cr{rr}}(\sigma_3)$ it holds that
\begin{equation}
\label{eq:Feq}
\min\bigl(\Cr{rr},{\rm adm}^-(\zeta(j),t)\bigr)=\min\bigl(\Cr{rr},|\{k\in g^{-1}(j): \zeta(k)\prec t\}|\bigr).
\end{equation}

Assume for contradiction that the construction gets stuck when trying to extend the domain of $g$ to some $i\in [\Cr{N}]$. Let $\tau=\zeta(i)$, let $t_1=\pi_{\Cr{rr}}(\tau)$, and let $t_2$ be the unique rank $\Cr{rr}$ local type such that $\varphi_{t_2}(f(x))\in \tau_1$.
By assumption, for every $\tau_2\in T_{\Cr{clean}}$ with $\tau_2\prec t_2$ it holds that ${\rm adm}^-(\tau_2,t_1)\leq \Cr{rr}$. Hence, by the $(\Cr{clean},\Cr{rr})$-restricted MTP, it holds that
$$
\sum_{\tau_1\prec t_1} {\rm adm}^+(\tau_1,t_2)\mu(\tau_1)=\sum_{\tau_2\prec t_2}{\rm adm}^-(\tau_2,t_1).
$$
Thus
$$
|\{i: \pi_{\Cr{rr}}(\zeta(i))=t_1\}|=\sum_{j}|\{k\in g^{-1}(j): \zeta(k)\prec t_1\}|,
$$
which contradicts the hypothesis that the construction gets stuck.

Now assume for contradiction that \eqref{eq:Feq} does not hold. Then there exists $t_1$ and $j_0$ such that 
$$
|\{k\in g^{-1}(j_0): \zeta(k)\prec t_1\}|<\min\bigl(\Cr{rr},{\rm adm}^-(\zeta(j_0), t_1)\bigr).$$
Let $t_2=\pi_{\Cr{rr}}(\zeta(j_0))$.
According to the construction of $g$, it holds for every $j$ such that $\zeta(j)\prec t_2$ that
$$|\{k\in g^{-1}(j): \zeta(k)\prec t_1\}|\leq \min\bigl(\Cr{rr},{\rm adm}^-(\zeta(j), t_1)\bigr).$$
Hence we have
\begin{align*}
\sum_{\tau_2\prec t_2}\min\bigl(\Cr{rr},{\rm adm}^-(\zeta(j),t)\bigr)\widehat\mu(\tau_2)&>\frac{1}{\Cr{N}}
\sum_{\tau_2\prec t_2}\sum_{j: \zeta(j)=\tau_2}|\{k\in g^{-1}(j): \zeta(k)\prec t_1\}|\\
&=\frac{1}{\Cr{N}}\sum_{\tau_1\prec t_1}\sum_{i: \zeta(i)=\tau_1}{\rm adm}^+(\zeta(i),t_2)\\
&=\sum_{\tau_1\prec t_1}{\rm adm}^+(\tau_1,t_2)\,\widehat\mu(\tau_1)
\end{align*}
which contradicts the $(\Cr{clean},\Cr{rr})$-restricted MTP. Thus \eqref{eq:Feq} holds.

The $\sigma_3$-mapping $\mathbf F_3$ has domain $[\Cr{N}]$. For every unary relation symbol $S\in \sigma_3$ we let
$S(\mathbf F_3)=\{i\in F_3:\ S(x)\in \zeta(i)\}$, and define
$f_{\mathbf F_3}=g$. 

Note that $\mathbf F$ has no cycle of length $\ell$ with $1\leq \ell\leq \Cr{cut}$: as $f(x)\wedge U_i(x)\rightarrow U_{(i+1) \bmod \Cr{cut}}(f(x))$ holds with probability $1$.
Hence, all the cycles have their length a multiple of $\Cr{cut}$.

That $\Type{\mathbf F}{\Cr{rr}}(v)=\pi_{\Cr{rr}}(\zeta(v))$ holds for every $v\in F_3$ then follows from Lemma~\ref{lem:lef}.
\end{proof}

As a consequence of Lemma~\ref{lem:ratmu} and Equation~\ref{eq:F3} it holds that
\begin{equation}
	\label{eq:distF3}
	\ldist{1}{\Cr{rr}}(\mathbf F_3,\mathbf L_3)<\Cr{epsmu}.
\end{equation}
\subsection{From $\mathbf F_3$ to $\mathbf F_2$: Rewiring short cycles}

We now let $\mathbf F_2=\mathsf I_2(\mathbf F_3)$.
Every  local formula $\phi\in{\rm FO}_1^{\rm local}(\sigma_2)$ with local rank at most $2\Cr{r}$ corresponds (for the $\mathsf I_2$ interpretation) to a local formula $\widehat\phi$ with local rank at most $2\Cr{r}(2\Cr{r}-1)<\Cr{rr}$. Then it holds that
$$
	|\langle\phi,\mathbf F_2\rangle-\langle\phi,\widetilde{\mathbf L}_2\rangle|=
	|\langle\widehat\phi,\mathbf F_3\rangle-\langle\widehat\phi,\mathbf L_3\rangle|.
$$
Thus
\begin{equation}
	\label{eq:distF2}
	\ldist{1}{2\Cr{r}}(\mathbf L_2,\mathbf F_2)\leq 
	\ldist{1}{2\Cr{r}}(\mathbf L_2,\widetilde{\mathbf L}_2)+
		\ldist{1}{2\Cr{r}}(\widetilde{\mathbf L}_2,\mathbf F_2)<\Cr{epsmu}.
\end{equation}
\subsection{The mapping $\mathbf E_1$: A finite model}
	
A {\em terminal} of $T_R$ is a type $\tau$ such that if $t'$ is such that ${\rm adm}^+(\tau,t')=1$ then
$\sum_{\tau'\prec t'}\mu(\tau')=0$. Importance of terminal types will be a consequence of  the following useful fact:
\begin{claim}
	Let $\tau_1$ be such that $\mu(\tau_1)>0$, and let $t_2$ be such that ${\rm adm}^+(\tau_1,t_2)=1$. 
	
	Then at least one of the following holds:
	\begin{enumerate}
		\item there exists $\tau_2\prec t_2$ such that $\mu(\tau_2)>0$;
		\item there exists $\tau_2\prec t_2$ such that  ${\rm adm}^-(t_1,\tau_2)>r$.
	\end{enumerate}
\end{claim}
\begin{proof}
Let $t_1$ be such that $\tau_1\prec t_1$.
	Assume that for every $\tau_2\prec t_2$ such that  ${\rm adm}^+(\tau_1,t_2)=1$ it holds  ${\rm adm}^-(t_1,\tau_2)\leq r$. Then, according to the FMTP, it holds
\begin{align*}
r\sum_{\vartheta_2\prec t_2}\mu(\vartheta_2)
&\geq
	\sum_{\vartheta_2\prec t_2}{\rm adm}^-(\vartheta_2,t_1)\mu(\vartheta_2)\\
	&=	\sum_{\vartheta_1\prec t_1}{\rm adm}^+(\vartheta_1,t_2)\mu(\vartheta_1)\\
	&\geq \mu(\tau_1)>0
\end{align*}
Thus there exists $\tau_2\prec t_2$ such that $\mu(\tau_2)>0$.
	\end{proof}

A type $\tau'$ is a {\em hub type} if there exists $\tau\prec t$ such that $\tau$ is a terminal and ${\rm adm}^-(\tau',t)>r$.
	Let $\tau_1,\dots,\tau_k$ be the terminal types of $\mathbf L$, and let $\tau_1',\dots,\tau_k'$ be associated hub types.
	
\begin{lemma}
\label{lem:model}
There exists a finite mapping $\mathbf M$ such that 
$\mathbf M\equiv_{\Cr{elem}}\mathbf L$, and such that there are
elements 
$$h_{1,1},\dots,h_{1,\Cr{away}},\dots,h_{k,1},\dots,h_{k,\Cr{away}}\in M,$$
 pairwise at distance at least $2^r$, such that $\Type{\mathbf M}{\Cr{rr}}(h_{i,j})=\tau_i'$.
\end{lemma}
\begin{proof}
	We consider the formula $\zeta$ with free variables
	$$x_{1,1},\dots,x_{1,\Cr{away}},\dots,x_{k,1},\dots,x_{k,\Cr{away}},$$ defined by
$$	\zeta:=
	 \Bigl(\bigwedge_{(i,j)\neq (i',j')}{\rm dist}(x_{i,j}, x_{i',j'})>2\Cr{r}\Bigr)\wedge\Bigl(\bigwedge_{1\leq i\leq k}\bigwedge_{1\leq j\leq \Cr{away}}\varphi_{\tau_i'}(x_{i,j})\Bigr)
$$
and the sentence 
$$
\theta:=(\exists x_{1,1},\dots,x_{1,\Cr{N}},\dots,x_{k,1},\dots,x_{k,\Cr{away}})\zeta
$$
The hub types can be chosen in such a way that $\theta$ is satisfied in $\mathbf L$. Indeed, for each $\tau'<t'$ such that ${\rm adm}(\tau,t')=1$ the connected component of any $v\in\phi_{\tau'}(\mathbf L)$ has measure $0$ (as $\mathbf L$ is residual) hence it is possible, for each terminal $\tau$ to choose $\tau'$ in such a way that there are in $\mathbf L$ uncountably many connected components with an element in $\phi_{\tau'}(\mathbf L)$.
\end{proof}

\subsection{From $\mathbf E_1$ and $\mathbf F_2$ to 
$\mathbf F_1$: Merging}

Let $S=\{v_1,\dots,v_k\}$ be the set of all terminal elements of $\mathbf F_2$, and let $\gamma(v_i)$ be the rank $r$-type corresponding to elements of $\mathbf L$ having type $\Upsilon(v_i)$ in $\mathbf L_1$.
 
\begin{figure}[h]
	\begin{center}
		\includegraphics[width=.7\textwidth]{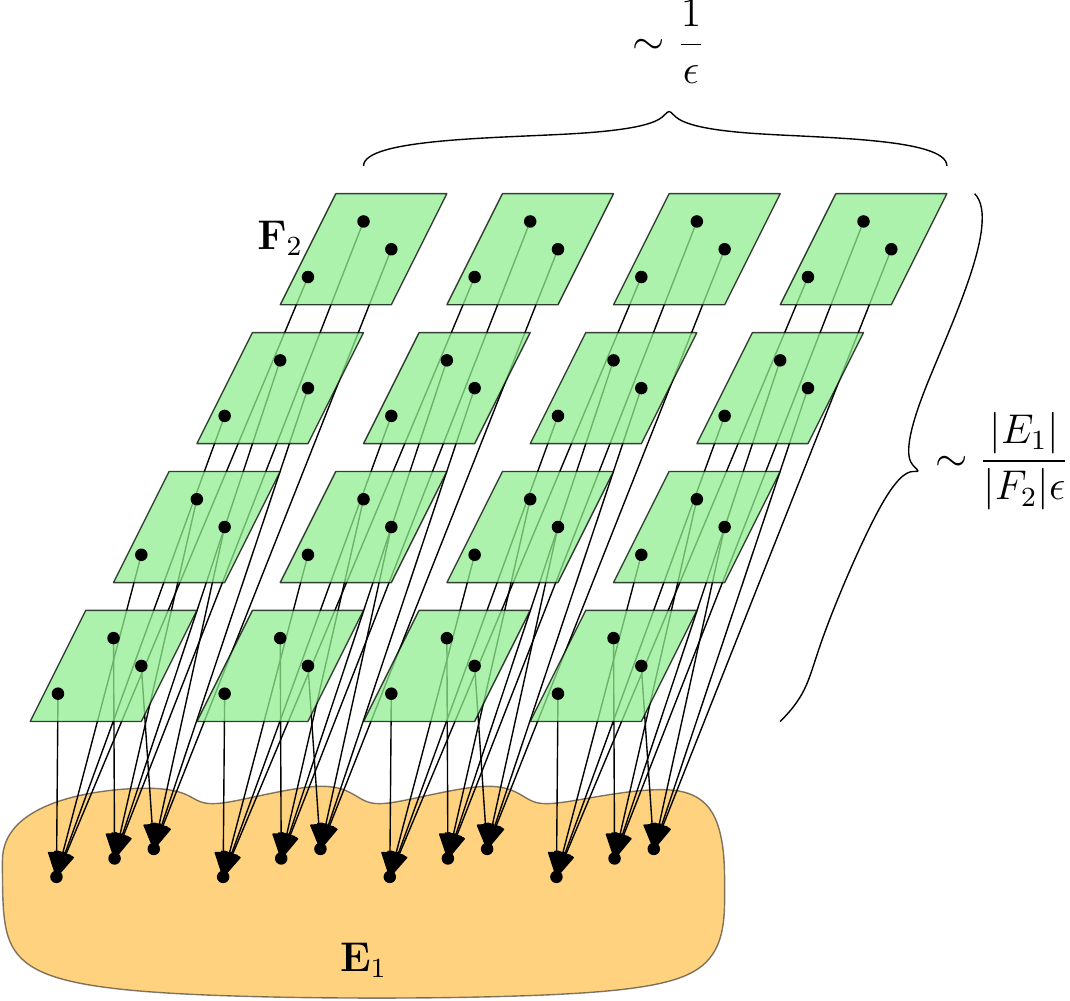}
	\end{center}
	\caption{Merging $\mathbf M$ with copies of $\mathbf F$}
	\label{fig:merge}
\end{figure}

 Let $\Cl[N]{close}=\lceil\frac{|E_1|}{|F_2|\Cr{res}}\rceil$, 
 and  let $F_1$ be the disjoint union of $E_1$ and 
$F_2\times[\Cr{close}]\times[\Cr{away}]$. If $v\in M$, we define $f_{\mathbf F_1}(v)=f_{\mathbf E_1}(v)$. Otherwise, if
$(v,i,j)\in F_2\times[\Cr{close}]\times[\Cr{away}]$ we define

$$
f_{\mathbf F_1}(v,i,j)=\begin{cases}
	(f_{\mathbf F_2}(v),i,j)&\text{if }v\notin S\\
	h_{a,i}&\text{if }v=v_a\in S
\end{cases}
$$
(See Fig.~\ref{fig:merge}.)

	We consider the finite mapping $\widetilde{\mathbf E}$ obtained from $\mathbf E_1$ as follows:
	For $1\leq i\leq k$ and $1\leq j\leq \Cr{close}$, and every 
	$z\in E_1$ such that $f_{\mathbf E_1}(z)=h_{i,j}$ and $\Type{\mathbf E_1}{\Cr{r}}(z)=\gamma(v_i)$, we mark $z$ by mark $A_{i,j}$ and let $f_{\widetilde{\mathbf E}}(z)=z$. For all other elements $z\in E_1$ we let $f_{\widetilde{\mathbf E}}(z)=f_{\mathbf E_1}(z)$. Moreover, each $h_{i,j}$ receives mark $B_{i,j}$.
	There is an easy basic quantifier-free interpretation $\mathsf I$ such that
	$\mathsf I(\widetilde{\mathbf E})=\mathbf E_1$.
	
	Now we consider the disjoint union $\widetilde{\mathbf F}$ of $\widetilde{\mathbf E}$ and $\Cr{away}\Cr{close}$ copies of $\mathbf F_2$, such that terminal $v_i$ in copy $(j,k)$ is marked $A_{i,j}$, and we let $\mathbf F_1=\mathsf I(\widetilde{\mathbf F})$.

\begin{lemma}
	The finite mappings $\mathbf E_1$ and $\mathbf F_1$ are $\Cr{r}$-equivalent.
\end{lemma}
\begin{proof}
	It is a direct consequence of Hanf's locality theorem that 
	$\widetilde{\mathbf F}$ is $\Cr{r}$-equivalent to $\widetilde{\mathbf E}$. It follows that $\mathbf F_1=\mathsf I(\widetilde{\mathbf F})$ is $\Cr{r}$-equivalent to ${\mathbf E_1}=\mathsf I(\widetilde{\mathbf E})$.
\end{proof}
\begin{lemma}
\label{lem:goodt}
	Each element $(v,i,j)$ in a copy of $\mathbf F_2$ in $\mathbf F_1$ is such that 
	$$\Type{\mathbf F_1}{\Cr{r}}(v,i,j)=\gamma(v).$$
\end{lemma}
\begin{proof}
	This follows easily from an Ehrefeucht-Fra{\"\i}ss{\'e} game.
\end{proof}
\begin{lemma}
$$
\langle[{\rm dist}(x_1,x_2)\leq 2\Cr{r},\mathbf F_1\rangle<\Cr{res}
$$
\end{lemma}
\begin{proof}
	Every ball of radius $2\Cr{r}$ contains less than $|E_1|+\Cr{close}|F_2|$ elements.
	Thus the probability $\langle[{\rm dist}(x_1,x_2)\leq 2\Cr{r},\mathbf F_1\rangle$ that two random elements in $\mathbf F_1$ are at distance at most $2\Cr{r}$ is less than $\frac{|E_1|+\Cr{close}|F_2|}{|E_1|+\Cr{away}\Cr{close}|F_2|}<\Cr{res}$.
\end{proof}

\begin{lemma}
It holds that
\begin{equation}
	\label{eq:distF1}
	\ldist{1}{\Cr{r}}(\mathbf F_1,\mathbf L_1)<\Cr{F1}.
\end{equation}
\end{lemma}
\begin{proof}
Let  $\phi\in{\rm FO}_1^{\rm local}$ be a formula with local rank at most $\Cr{r}$. 	Let $\psi=\bigvee_{t\ni\phi}R_t$, where the disjunction is over rank $\Cr{r}$-local types.
Then $\langle\phi,\mathbf L_1\rangle=\rangle\psi,\mathbf L_2\rangle$.
According to Lemma~\ref{lem:goodt} it holds that
$$
|\langle\phi,\mathbf F_1\rangle-\langle\psi,\mathbf F_2\rangle|\leq \frac{|E_1|}{|F_1|}\leq \frac{1}{1+\Cr{away}\Cr{close}\frac{|F_2|}{|E_1|}}\leq\frac{1}{1+\frac{2}{\Cr{res}^2}}<\frac{\Cr{res}^2}{2}.
$$
Thus
\begin{align*}
|\langle\phi,\mathbf F_1\rangle-\langle\phi,\mathbf L_1\rangle|&\leq
|\langle\phi,\mathbf F_1\rangle-\langle\psi,\mathbf F_2\rangle|+
|\langle\psi,\mathbf F_2\rangle-\langle\psi,\mathbf L_2\rangle|\\
&<\frac{\Cr{res}^2}{2}+\Cr{epsmu}<\Cr{F1}.
\end{align*}
\end{proof}

\subsection{From $\mathbf F_1$ to $\mathbf F$: approximation of the original mapping}

At this stage, we have constructed a finite mapping $\mathbf F_1$
such that $\mathbf L_1 \equiv_{\Cr{r}} \mathbf F_1$ and $|\langle \psi, \mathbf L_1 \rangle - \langle \psi, \mathbf F_1 \rangle| < \Cr{F1}$ for every $\psi \in {\rm FO}^{\rm local}_1$ with rank at most $\Cr{r}$.

Let $\mathbf F=\mathsf I(\mathbf F_1)$, where $\mathsf I_1$ is the interpretation defined in Section~\ref{sec:L1}.
The following lemma ends the proof of Theorem~\ref{thm:invfomap}.
\begin{lemma}
For every formula $\phi$ with $p$ free variables and rank at most $r$ it holds that
 $$|\langle \varphi, \mathbf L \rangle - \langle \varphi, \mathbf F \rangle| < \epsilon.$$
\end{lemma}
\begin{proof} 
Let $\phi$ be a local formula with at most $p$ free variables and rank at most $\Cr{r}$.

As $\mathbf L_1$ is $\Cr{res}$-residual, according to Lemma~\ref{lem:res}, it holds that

\begin{align*}
\ldist{p}{\Cr{r}}(\mathbf L_1,\mathbf F_1)&\leq 2p\ldist{1}{\Cr{r}}(\mathbf L_1,\mathbf F_1)+\binom{p}{2}(\langle\delta_{2\Cr{r}},\mathbf L_1\rangle+\langle\delta_{2\Cr{r}},\mathbf F_1\rangle)\\
&<2p\Cr{F1}+\binom{p}{2}\Cr{res}<\Cr{eps}.
\end{align*}

We deduce from
Lemma~\ref{lem:locred} and the definitions of $\Cr{r}$ and $\Cr{eps}$ that 
$|\langle \widetilde\varphi, \mathbf L_1 \rangle - \langle \widetilde\varphi, \mathbf F_1 \rangle| < \epsilon$
holds true for every first-order formula $\widetilde\varphi$ with at most $p$ free variables and rank at most $r$.

Let $\varphi$ be a first-order formula with at most $p$ free variables and rank at most $r$. Then there exists a formula
$\widetilde\varphi$ with at most $p$ free variables and rank at most $r$ such that $\langle \widetilde\varphi,\mathbf L_1\rangle=\langle \varphi,\mathbf L\rangle$ and $\langle \widetilde\varphi,\mathbf F_1\rangle=\langle \varphi,\mathbf F\rangle$. Hence
$|\langle \varphi, \mathbf L \rangle - \langle \varphi, \mathbf F \rangle| < \epsilon$.
\end{proof}

This ends the last reduction step in the proof of Theorem~\ref{thm:invfomap}. As explained above (see Fig.~\ref{fig:stratFO} and comments preceding it) this finishes the proof of Theorem~\ref{thm:invfomap}.

\section{Local approximation}

The aim of this section is to prove Theorem~\ref{thm:invlocmap} by following steps similar to those we followed to prove Theorem~\ref{thm:invfomap}.

The first main difference is that we cannot use general first-order interpretations, but only local interpretations. Thus we cannot follow the first reduction step to reduce to the $\epsilon$-residual case. Instead, we shall prove that every connected mapping modeling is close (for the topology of local convergence) to a connected mapping modeling with the finite model property, for which Theorem~\ref{thm:invfomap} applies.
The strategy will be to consider first the connected components of $\mathbf L$ with non-negligible measures, and then the remaining components of the mapping modeling.

For $\epsilon$-residual mapping modelings, we can follow the proof of Theorem~\ref{thm:invfomap} until Step~\ref{step:merge}. In this step, the model $M$ will be replaced by the union of models of the hub local types.

\subsection{Connected mapping modelings}

Let $\mathbf L$ be a connected mapping modeling. We define a directed graph modeling $\widehat{\mathbf L}$ with countably many marks $M$ and $N$ as follows:
\begin{itemize}
	\item The domain of $\widehat{\mathbf L}$ is $L$, with same probability measure;
	\item if $Z(\mathbf L)\neq\emptyset$, we arbitrarily mark a vertex $v\in Z(\mathbf L)\neq\emptyset$ with mark $M$ and its image $f(v)$ with mark $N$;
	\item the arcs of $\widehat{\mathbf L}$ are the pairs $(v,f(v))$ for which $v$ is not marked by $M$.
\end{itemize}

The following lemma is much stronger than what we need. It would be sufficient to say that for some $d$ the ball of radius $d$ around the root has measure at least $1-\epsilon$.
Now the idea is that the ball $B$ of radius $d+r$ around the root of $\mathbf L$ not only has measure close to $1$, but also has the property that less than $\epsilon$ measure of the elements have different rank $r$ local type in $\mathbf L$ and $\mathbf L\mid B$. Now $\mathbf L\mid B$ has finite height hence enjoys the finite model property. An FO-approximation of $\mathbf L\mid B$ is then an ${\rm FO}^{\rm local}$-approximation of $\mathbf L$.

\begin{lemma}
	\label{lem:cut}
	Let $\mathbf L$ be a connected mapping modeling with atomless measure $\nu_{\mathbf L}$ that satisfies the MTP, and let $\epsilon>0$ be a positive real.
	
	Then for every $r\in L$ there exists $d\in\bbbn$
	such that the the subset $A\subseteq L$, defined as the union of the vertex sets of all the (undirected) paths of length at least $d+1$ in $\widehat{\mathbf L}$ with endpoint $r$, has measure at most $\epsilon$.
\end{lemma}
\begin{proof}
	There exists an even integer $d$ such that the ball $B_{d/2}(\widehat{\mathbf L},r)$ has measure at least $(1-\epsilon/2)$. 
	 For $0\leq i\leq d$, let
	$S_i$ be the set of all vertices of $A$ at distance exactly $i$ from $r$.
	According to the MTP (and uniqueness of paths from a vertex $v$ to $r$), and as $\nu_{\mathbf L}$ is atomless,
	it holds that
	$$
	0=\nu_{\mathbf L}(\{r\})\leq\nu_{\mathbf L}(S_{1})\leq\dots\leq\nu_{\mathbf L}(S_{d}).
	$$
	Thus it holds that
	\begin{align*}
	\nu_{\mathbf L}\Bigl(\bigcup_{i=0}^{d/2}S_i\bigr)&\leq 
		\nu_{\mathbf L}\Bigl(\hspace{-3mm}\bigcup_{{i=d/2+1}}^{d}\hspace{-3mm}S_i\bigr).
		\intertext{That is:}
	\nu_{\mathbf L}(A\cap B_{d/2}(\widehat{\mathbf L},r))&\leq 
	\nu_{\mathbf L}(\bigcup_{i=d/2}^d S_{i})\\
	&\leq \nu_{\mathbf L}(L\setminus B_{d/2}(\widehat{\mathbf L},r)).	\intertext{Thus}
	\nu_{\mathbf L}(A)&\leq \nu_{\mathbf L}(A\cap B_{d/2}(\widehat{\mathbf L},r))+\nu_{\mathbf L}(A\setminus B_{d/2}(\widehat{\mathbf L},r))\\
	&\leq 2\nu_{\mathbf L}(L\setminus B_{d/2}(\widehat{\mathbf L},r))\\
	&<\epsilon
	\end{align*}
\end{proof}

\begin{definition}
	Let $\mathbf L$ be a colored mapping modeling with finite height and let $r\in\bbbn$.
	We define the {\em standard $r$-approximation} $\widehat{\mathbf L}$ of $\mathbf L$ as follows:
	
	Let $C=Z(\mathbf L)$ and $C_i=Z_i(\mathbf L)$. For $x\in L$ let $h(x)$ be the minimum non-negative integer $k$ such that $f_{\mathbf L}^k(x)\in C$. Note that $0\leq h(x)\leq {\rm height}(\mathbf L)$. Let $p=\max_{x\in C\setminus C_1} h(x)$. We iteratively define sets $X_i$ for $i=p,\dots,1$, together with an equivalence relation $\sim_i$ on $h^{-1}(i)\cap f_{\mathbf L}(X_{i+1})$ (if $i<p$). We start with $i=p$ and define $\sim_{p}$ on $L$ by $x\sim_{p+1} y$ if $x$ and $y$ have the same color. For every 
		$y\in h^{-1}(i-1)$ we choose an inclusion maximal subset $I(y)$ of $f_{\mathbf L}^{-1}(y)$ containing no $r+1$ $\sim_i$-equivalent vertices. Then we define $X_i=\bigcup_{y\in h^{-1}(i-1)}I(y)$, and we define the equivalence relation 
	$\sim_{i-1}$ on $h^{-1}(i-1)\cap f_{\mathbf L}(X_{i})$ by $y_1\sim_{i-1}y_2$ if for every $z\in f_{\mathbf L}^{-1}(y_1)\cup f_{\mathbf L}^{-1}(y_2)$ it holds that
	$$
	|\{x_1\in f_{\mathbf L}^{-1}(y_1): x_1\sim_i z\}|=
		|\{x_2\in f_{\mathbf L}^{-1}(y_2): x_2\sim_i z\}|.
	$$
	We now consider the restriction $g$ of $f_{\mathbf L}$ to $C\cup\bigcup_{i=1}^p X_i$. Note that all the connected components have their size bounded by some fixed function of $c$ and $p$. We consider an inclusion maximal union $\widehat{\mathbf L}$ of connected components of $g$ containing no $r+1$ isomorphic connected components. The mapping $\widehat{\mathbf L}$ is then the restriction of $\mathbf L$ to $\widehat{\mathbf L}$. Note that $\widehat{\mathbf L}$ has its size bounded by some fixed function of $c$ and $p$.
 \end{definition}
 
 An alternate construction can be used, which is parametrized by a pair $(r,R)$ of integers with $r\leq R$.
 The idea is as follows: we start from the standard $R$-approximation and then reduce every set of at least $k>r$ equivalent sons to $r$ if either $k<R$ or some descendent of one of these sons as $R$ equivalent sons.
 Then, according to MTP, the measure of the types of the vertices obtained by removing any $R$ equivalent siblings and their descendants is at most $F(r,t)/R$. So one should require $R> F(r,t)/\epsilon$.

\begin{lemma}
\label{lem:fmpheight}
	Every mapping $\mathbf L$ with finite height 
	is $r$-equivalent to its standard $r$-approximation, hence
	has the finite model property.
\end{lemma}
\begin{proof}
An easy strategy for the $r$-round Ehrenfeucht-Fra{\"\i}ss{\'e} game shows that $\mathbf L$ is $r$-equivalent to its standard $r$-approximation.
	\end{proof}

\subsection{Merging with hub local type models}
To each rank $r$ hub local type $\tau_i'$ we associate a finite rooted mapping $(\mathbf M_i,h_i)$ such that $\Type{\mathbf M_i}{r}(h_i)=\tau_i'$. Let $\mathbf M$ be the disjoint union of the $\mathbf M_i$. We proceed to the merge of $\mathbf M$ with copies of $\mathbf F$ as in Step~\ref{step:merge} of the proof of Theorem~\ref{thm:invfomap}.

\begin{figure}[h]
	\begin{center}
		\includegraphics[width=.7\textwidth]{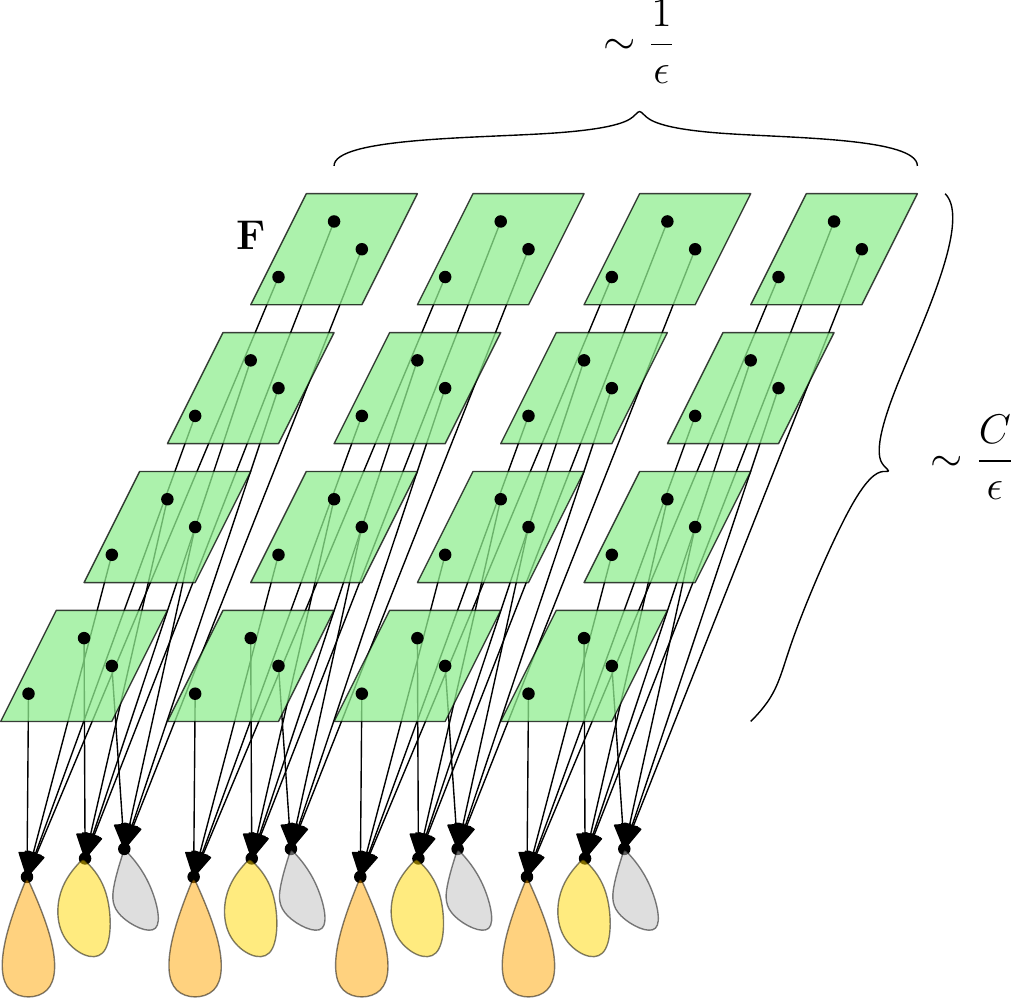}
	\end{center}
	\caption{Merging small models with many copies of $\mathbf F$.}
	\label{fig:merge2}
\end{figure}

\section{Concluding Remarks}
In this paper we considered the approximation problem for mapping modelings. It would be interesting to consider the approximation problem where we have 
 only the probability measure $\mu$ corresponding to the satisfaction probability of first-order formulas.
 
 In such a setting, we shall consider probability measures $\mu$ on $S(\mathcal L_{\rm FO})$ that are invariant under the natural action of the infinite permutation group $S_\omega$ (acting by permuting the free variables in the formulas), whose support projects on a single point ${\rm Th}(\mu)$ of $S(\mathcal L_{{\rm FO}_0})$. The analogs of the property we required for modeling mappings are as follows:
  The condition for the modeling to be atomless corresponds to the property that the $\mu$-measure of the clopen subset $K(x_1=k_2)$ of $S(\mathcal L_{\rm FO})$ dual to the formula $x_1=x_2$ is zero.
The  finite model property of the modeling corresponds to the property that every sentence in ${\rm Th}(\mu)$ has a finite model.
The finitary mass transport principle for the modeling corresponds to the following property of $\mu$: for every formulas $\phi,\psi\in {\rm FO}_1$ such that $\psi(x)$ entails that there exist exactly (resp. strictly more than) $k$ elements $y_1,\dots,y_k$ such that $\phi(y_i)\wedge f(y_i)=x$ we have $\mu(K(\phi))=k\mu(K(\psi))$ (resp. $\mu(K(\phi))>k\mu(K(\psi))$).

Admittedly the proofs presented in this paper are technical and complex. 
In a way this was expected as 
approximating modeling structures with two mappings seem to be fully out of reach. Indeed, the existence of a finite (local) approximation for modelings consisting into two (bjiective) mappings $f$ and $g$ with $f^2=g^3={\rm Id}$ satisfying the FMTP is equivalent to the general Aldous-Lyons conjecture.

An interesting question is to solve the inverse problem for acyclic modelings (the modeling equivalent of {\em treeings}). This problem has been solved in the bounded diameter case \cite{limit1} by a complicated analysis, and in the bounded degree case by \cite{Elek2010}.
However the problem for general acyclic modelings remain open.

A way to make the problem simpler is to assume that the acyclic modeling looks like a directed rooted tree. This is the motivation of the following problem stated in \cite{modeling}:
 if a tree modeling is oriented in such a way that the root is a sink and non-roots have outdegree one and if any finite subset of the complete theory of the modeling has a connected finite model, is it true that the modeling is the FO-limit of a sequence of finite rooted trees?

Finally, we would like to mention that random mappings are not FO-convergent, as they do not satisfy a 0-1 law (the expected number of cycles of length $r$ tend to $1/r$ \cite{flajolet1989random}). However it might be possible that random mappings are ${\rm FO}^{\rm local}$-convergent.
\providecommand{\bysame}{\leavevmode\hbox to3em{\hrulefill}\thinspace}
\providecommand{\MR}{\relax\ifhmode\unskip\space\fi MR }
\providecommand{\MRhref}[2]{%
  \href{http://www.ams.org/mathscinet-getitem?mr=#1}{#2}
}
\providecommand{\href}[2]{#2}

\end{document}